\documentclass[a4paper,12pt]{amsart}
\usepackage{amsfonts}

\usepackage{amssymb}
\usepackage{ifthen}
\usepackage{graphicx}
\usepackage{tablists}
\nonstopmode \numberwithin{equation}{section}
\usepackage[usenames]{color}
\usepackage{enumerate}
\newtheorem{thm}[equation]{Theorem}
\newtheorem{lem}[equation]{Lemma}
\newtheorem{cor}[equation]{Corollary}

\newtheorem{cl}{Claim}[section]
\newtheorem{ca}{Case}[section]
\newtheorem{sca}{Subcase}[section]
\newtheorem{scl}[section]{Subclaim}
\newtheorem{conj}[equation]{Conjecture}

\theoremstyle{definition}
\newtheorem{defn}[equation]{Definition}
\newtheorem{prob}[equation]{Problem}
\newtheorem{op}[equation]{Open Problem}
\newtheorem{ques}[equation]{Question}
\newtheorem{rem}[equation]{Remark}
\newtheorem{exam}[equation]{Example}

\newcounter {own}
\def\theown {\thesection       .\arabic{own}}

\newenvironment{pf}[1][]{%
 \vskip 3mm
 \noindent
 \ifthenelse{\equal{#1}{}}%
  {{\slshape Proof. }}%
  {{\slshape #1.} }%
 }%
{\qed\bigskip}

\newcommand{\IR}{{\mathbb R}}

\newcommand{\dist}{{\operatorname{dist}}}

\def\be{\begin{equation}}
\def\ee{\end{equation}}

\newcommand{\ben}{\begin{enumerate}}
\newcommand{\een}{\end{enumerate}}

\newcommand{\blem}{\begin{lem}}
\newcommand{\elem}{\end{lem}}
\newcommand{\bthm}{\begin{thm}}
\newcommand{\ethm}{\end{thm}}
\newcommand{\bcor}{\begin{cor}}
\newcommand{\ecor}{\end{cor}}
\newcommand{\beg}{\begin{exam}}
\newcommand{\eeg}{\end{exam}}
\newcommand{\begs}{\begin{examples}}
\newcommand{\eegs}{\end{examples}}
\newcommand{\bdefe}{\begin{defn}}
\newcommand{\edefe}{\end{defn}}
\newcommand{\bprob}{\begin{prob}}
\newcommand{\eprob}{\end{prob}}
\newcommand{\bques}{\begin{ques}}
\newcommand{\eques}{\end{ques}}
\newcommand{\bei}{\begin{itemize}}
\newcommand{\eei}{\end{itemize}}
\newcommand{\bcon}{\begin{conj}}
\newcommand{\econ}{\end{conj}}
\newcommand{\bop}{\begin{op}}
\newcommand{\eop}{\end{op}}

\newcommand{\bca}{\begin{ca}}
\newcommand{\eca}{\end{ca}}
\newcommand{\bsca}{\begin{sca}}
\newcommand{\esca}{\end{sca}}

\newcommand{\bcl}{\begin{cl}}
\newcommand{\ecl}{\end{cl}}

\newcommand{\bscl}{\begin{scl}}
\newcommand{\escl}{\end{scl}}

\newcommand{\bcons}{\begin{conjs}}
\newcommand{\econs}{\end{conjs}}
\newcommand{\bprop}{\begin{propo}}
\newcommand{\eprop}{\end{propo}}
\newcommand{\br}{\begin{rem}}
\newcommand{\er}{\end{rem}}
\newcommand{\brs}{\begin{rems}}
\newcommand{\ers}{\end{rems}}
\newcommand{\bo}{\begin{obser}}
\newcommand{\eo}{\end{obser}}
\newcommand{\bos}{\begin{obsers}}
\newcommand{\eos}{\end{obsers}}
\newcommand{\bpf}{\begin{pf}}
\newcommand{\epf}{\end{pf}}
\newcommand{\ba}{\begin{array}}
\newcommand{\ea}{\end{array}}
\newcommand{\beq}{\begin{eqnarray}}
\newcommand{\beqq}{\begin{eqnarray*}}
\newcommand{\eeq}{\end{eqnarray}}
\newcommand{\eeqq}{\end{eqnarray*}}

\newcounter{minutes}
\divide\time by 60
\newcounter{hours}
\multiply\time by 60 \addtocounter{minutes}{-\time}

\begin{document}

\bibliographystyle{amsplain}
\title{Busemann functions and uniformization of Gromov hyperbolic spaces}

\author{Qingshan Zhou}
\address{Qingshan Zhou, School of Mathematics and Big Data, Foshan University,  Foshan, Guangdong 528000, People's Republic
of China} \email{qszhou1989@163.com; q476308142@qq.com}

\author{Saminathan Ponnusamy}
\address{Saminathan Ponnusamy, Department of Mathematics, Indian Institute of Technology Madras, Chennai 600036,
India.}
\address{Lomonosov Moscow State University,
Moscow Center of Fundamental and Applied Mathematics, Moscow, Russia.
}
\email{samy@iitm.ac.in}

\author{Antti Rasila}
\address{Antti Rasila, Department of Mathematics with Computer Science, Guangdong Technion,  241 Daxue Road, Shantou, Guangdong 515063, People's Republic of China and Department of Mathematics, Technion -- Israel Institute of Technology, Haifa 32000, Israel} \email{antti.rasila@gtiit.edu.cn; antti.rasila@iki.fi}

\def\thefootnote{}
\footnotetext{ \texttt{\tiny File:~\jobname .tex,
          printed: \number\year-\number\month-\number\day,
          \thehours.\ifnum\theminutes<10{0}\fi\theminutes}
} \makeatletter\def\thefootnote{\@arabic\c@footnote}\makeatother

\date{}


\begin{abstract}
Uniformization theory of Gromov hypebolic spaces investigated by Bonk, Heinonen and Koskela, generalizes the case where a classical Poincar\'e ball type model is used as the starting point.
In this paper, we develop this approach in the case where the underlying domain is unbounded, corresponding to the classical Poincar\'e half-space model.
More precisely, we study conformal densities via Busemann functions on Gromov hyperbolic spaces and prove that the deformed spaces are unbounded uniform spaces. Furthermore, we show that there is a one-to-one correspondence between the bilipschitz classes of proper geodesic Gromov hyperbolic spaces that are roughly starlike with respect to a point on Gromov boundary and the quasisimilarity classes of unbounded locally compact uniform spaces. Our result can be understood as an unbounded counterpart of the main result of Bonk, Heinonen, and Koskela in ``Uniformizing Gromov Hyperbolic Spaces", Ast\'erisque 270 (2001).
\end{abstract}

\thanks{Qingshan Zhou was supported by NNSF of China (No. 11901090), by Guangdong Basic and Applied Basic Research Foundation (No. 2022A1515012441). Antti Rasila was supported by NNSF of China (No. 11971124) and NSF of Guangdong Province (No. 2021A1515010326).}

\maketitle{} \pagestyle{myheadings} \markboth{Zhou et al.}{Busemann functions and uniformization of Gromov hyperbolic spaces}

\section{Introduction and main results}\label{sec-1}

There are several models for the hyperbolic spaces $\mathbb{H}^n$, such as the Klein model, the Poincar\'{e} ball model, and the half-space model, see for example \cite[Chapter I.6]{BrHa}.
Study of these models provides characterizations of properties of these spaces such as geodesics, hyperplanes, and isometries. In \cite{Gr87}, Gromov showed that fundamental characteristics of $\mathbb{H}^n$ can be obtained by using a simple condition for quadruples of points, and defined a class of metric spaces known as Gromov hyperbolic spaces. They form a large and much studied class of metric spaces, which plays an important role in, for example, geometric group theory \cite{BrHa,Gr87}, analysis on metric spaces \cite{LS,ZLR}, and geometric function theory \cite{HSX,ZR}.

In \cite{BHK}, Bonk, et al. established the following Poincar\'{e} ball model for proper geodesic Gromov hyperbolic spaces.

\begin{thm}\label{t-1}$($\cite[Proposition 4.5]{BHK}$)$ The conformal deformations of $(X,d_\epsilon)$ of a proper geodesic Gromov hyperbolic space $(X,d)$ are bounded locally compact $A$-uniform spaces for some $\epsilon>0$.
\end{thm}

The notion of uniform domains was introduced by Martio and Sarvas in \cite{MS78}. This class of domains is a suitable generalization of the unit ball and quasidisks in the study of quasiconformal mapping theory. In the case of planar simply connected domains, several classical results can be generalized to higher dimensional Euclidean spaces of  quasiconformal mappings in uniform domains, see \cite{GS,MV11,MS78,ZP}.

Recently, Butler \cite{Bu1} obtained an unbounded counterpart for Gromov hyperbolic spaces conceptually similar to the Poincar\'{e} half-space model on $\mathbb{H}^n$, see \cite{Bu2} for
more applications on CAT($-1$) spaces and function spaces.
More precisely, it was shown in \cite[Theorem 1.3]{Bu1} that the conformal deformations of complete geodesic roughly starlike Gromov hyperbolic spaces are unbounded uniform spaces. For the definition of rough starlikeness, see Subsection \ref{sub-rs}.

In view of these considerations, the motivation of our study is to prove a result similar to Theorem \ref{t-1}, and to obtain uniformization of proper geodesic Gromov hyperbolic spaces. By introducing new conformal densities via Busemann functions on Gromov hyperbolic spaces, we show that the deformed spaces are unbounded uniform spaces.

\begin{thm}\label{z3} Suppose that $(X,d)$ is a proper geodesic $\delta$-hyperbolic space and the Gromov boundary $\partial_\infty X$ of $X$ contains at least two points. Then the conformal deformations $X_\kappa=(X,d_\kappa)$ induced by the densities \eqref{b-0} are unbounded locally compact $A$-uniform spaces with a constant $A=A(\delta)$ for all $0<\kappa\leq \kappa_0(\delta)$.
\end{thm}

\br Note that Theorem \ref{z3} is not covered in the work of Butler \cite[Theorem 1.3]{Bu1} because we do not assume the space $X$ to be roughly starlike. Not all geodesic Gromov hyperbolic spaces are roughly starlike. For example, let $\mathbb{H}^2$ be the hyperbolic plane and $o\in\mathbb{H}^2$. Attach a sequence of segments $\{[0, n]\}_{n\in\mathbb{N}}$ to $o$ by identifying $0$ with $o$. Then the resulting space $Y$, endowed with the induced length metric, is not roughly starlike.

We remark that this paper is closely related to an unpublished manuscript \cite{Z20} which is available in  arXiv:2008.01399. The present version is indeed a reformulated and improved version of this unpublished manuscript. It was also pointed out in \cite{Bu1} that the two works  \cite{Bu1} and \cite{Z20} were independently developed and essentially at the same time.
\er

For basic properties of Busemann functions we refer to \cite{BrHa,BuSc} and the references therein. The class of Busemann functions is an important tool in many research areas.
For example, Zhou, et al. \cite{ZPG} recently demonstrated that there is a quasisymmetric homeomorphism between the Euclidean boundary of an unbounded domain $G\subset \mathbb{R}^n$ and the punctured Gromov boundary of the space $(G,\widetilde{j}_G)$ equipped with a parabolic visual metric by using Busemann functions.

It is well known that the conformal mappings of the unit disk and the half plane onto itself are hyperbolic isometries and vice versa. In \cite{BHK}, Bonk, et al. generalized this result in
the metric spaces setting as follows.

\begin{thm}\label{BHK-thm}$($\cite[Theorem 1.1]{BHK}$)$ There is a one-to-one correspondence between the bilipschitz classes of proper, geodesic, and roughly starlike Gromov hyperbolic spaces and the quasisimilarity classes of bounded locally compact uniform spaces.
\end{thm}

This uniformization procedure plays an important role in the analysis of metric spaces \cite{BHK,GS,GPZ,LS,RSZ,ZR,ZR2} and the references therein. It allows us to replace Gromov hyperbolic spaces with uniform spaces where the geometry is more easier to understand. For example, Zhou and Rasila \cite{ZR} studied Teichm\"{u}ller's displacement problem on Gromov hyperbolic domains by using Theorem \ref{BHK-thm}.
Recently, Rogovin, et al. \cite{RSZ} obtained new characterizations of the Gehring-Hayman theorem, and used it to determine the critical exponents for the deformed spaces to be uniform spaces in the case of the hyperbolic spaces $\mathbb{H}^n$, the model spaces $\mathbb{M}^{\kappa}_n$ of the sectional curvature $\kappa<0$, and hyperbolic fillings.

As the second motivation, we consider whether there is a one-to-one correspondence between the bilipschitz classes of Gromov hyperbolic spaces and the quasisimilarity classes of unbounded uniform spaces. Applying Theorem \ref{z3}, we obtain the following main result:

\begin{thm}\label{main thm-1} There is a one-to-one correspondence between the bilipschitz classes of proper geodesic Gromov hyperbolic spaces that are roughly starlike with respect to points on the Gromov boundary and the quasisimilarity classes of unbounded locally compact uniform spaces.
\end{thm}

The paper is organized as follows. In Section \ref{sec-2} we recall some definitions.
Section~\ref{sec-3} focuses on Gromov hyperbolic spaces and we prove a few auxiliary lemmas. In Section~\ref{sec-4}, we introduce a class of metric densities by using Busemann functions, and then prove Theorem \ref{z3}. The proof of Theorem \ref{main thm-1} is given in Section \ref{sec-5}.

\section{Preliminaries}\label{sec-2}

\subsection{Notation} Following the notation of \cite{BuSc}, for $t_1,t_2,C\in \mathbb{R}$ with $C\geq 0$, we write $t_1\doteq t_2$ up to an error $\leq C$ or $t_1\doteq_C t_2$ instead of $|t_1-t_2|\leq C.$
For any two given sequences $\{t_i\},\{s_i\}\subset \mathbb{R}$, we write
$\{t_i\}_i \doteq \{s_i\}_i$ up to an error $\leq C$ or $\{t_i\}_i \doteq_C \{s_i\}_i$
if $\limsup\limits_{i\to \infty} |t_i-s_i|\leq C$. Similarly, we write
$\{t_i\}_i\leq \{s_i\}_i+C$ if
$\limsup\limits_{i\to \infty} t_i\leq \liminf\limits_{i\to \infty} s_i+C$,
and
$\{t_i\}_i\geq \{s_i\}_i+C$ if
$\liminf\limits_{i\to \infty} t_i\geq \limsup\limits_{i\to \infty} s_i+C.$

Let $A,B>0$ and $M\geq 1$. We use the notation $A\asymp B$ up to multiplicative error $\leq M$ or $A\asymp_M B$ instead of $1/M\leq A/B\leq M$.
For given two real numbers $s,t$, we set $s\wedge t=\min\{s,t\}$ and $s\vee t=\max\{s,t\}.$

\subsection{Metric geometry}
Let $(X,d)$ be a metric space. We say that $X$ is {\it incomplete} if its boundary $\partial X= \overline{X}\setminus X\not=\emptyset$, where $\overline{X}$ denotes the completion of $X$. Let $B(x,r)$  and $\overline{B}(x,r)$ be, respectively, the open ball and closed ball of radius $r$ centered at the point $x$ in $X$. The space $X$ is called {\it proper} if its closed balls are compact.

A {\it curve} in $X$ means a continuous map $\gamma:\, I\to X$ from an interval $I\subset \IR$ to $X$. If $\gamma$ is an embedding of $I$, then it is called an {\it arc}. We also denote the image set $\gamma(I)$ of $\gamma$ by  simply $\gamma$  itself. The {\it length} $\ell(\gamma)$ of $\gamma$ with respect to the metric $d$ is defined in an obvious way. Here the parametric interval $I$ is allowed to be open or half-open. We also denote the subarc of $\gamma$ by $\gamma[x,y]$ with endpoints $x$ and $y$ in $\gamma$.  Moreover, $X$ is called {\it rectifiably connected} if every pair of points in $X$ can be joined with a curve $\gamma$ in $X$ with $\ell(\gamma)<\infty$. The space $X$ is called {\it $A$-quasiconvex} if each pair of  points $x,y\in X$ can be joined by an $A$-{\it quasiconvex} curve $\gamma$, that is, $\ell(\gamma)\leq A d(x,y).$

Let $\gamma:I\to X$ be a curve. We say that $\gamma$ is a {\it geodesic arc} joining $x$ to $y$ in $X$ provided $\gamma(0)=x$, $\gamma(l)=y$, and $d(\gamma(t),\gamma(t'))=|t-t'|$ for all $t,t'\in I$. If $I=[0,\infty)$, then $\gamma$ is called a {\it geodesic ray}. If $I=\mathbb{R}$,  then $\gamma$ is called a {\it geodesic line}. A metric space $X$ is said to be {\it geodesic} if every pair of points can be joined by a geodesic arc. Denote by $[x,y]$ any geodesic between $x$ and $y$ in $X$.

Suppose that $(\Omega, d)$ is an incomplete, locally compact, and rectifiably connected metric space and that the identity map $(\Omega,d)\to (\Omega,\ell)$ is continuous, where $\ell$ is the intrinsic length metric of $\Omega$ with respect to $d$. As in \cite{BHK}, the {\it quasihyperbolic metric} $k$ in $\Omega$ is defined by
$$k(x, y)=\inf_{\gamma} \Bigg\{\int_{\gamma} \frac{1}{d(z)}\,|dz|\Bigg\},$$
where the infimum is taken over all rectifiable curves $\gamma$ in $\Omega$ with endpoints $x$ and $y$, $d(z):=\dist(z,\partial \Omega)$, and $|dz|$ denotes the arc-length element with respect to the metric $d$. We remark that $(\Omega,k)$ is proper and geodesic (cf. \cite[Proposition 2.8]{BHK}).

Let $A\geq 1$. An arc $\gamma$ connecting $x$ and $y$ in $\Omega$ is called {\it $A$-uniform}, if it is $A$-quasiconvex and for all $z\in \gamma$, $\ell(\gamma[x,z])\wedge \ell(\gamma[z,y])\leq A d(z).$
The space $\Omega$ is called $A$-{\it uniform} if each pair of points $x,y\in \Omega$ can be joined by an $A$-uniform arc.

  \subsection{Quasi-isometric and quasisimilarity maps}
Let $f:$ $(X,d)\to (X',d')$ be a map between metric spaces $X$ and $X'$, and let  $\lambda\geq 1$ and $C\geq 0$. We say that $f$ is $(\lambda, C)$-{\it quasi-isometric} if for all $x,y\in X$,
$$\lambda^{-1} d(x,y)-C\leq d'(f(x),f(y))\leq \lambda d(x,y)+C.$$
If in addition, every point $x'\in X'$ has distance at most $C$ from the set $f(X)$, then $f$ is called a {\it $(\lambda, C)$-quasi-isometry}. Moreover, if $f$ is a homeomorphism and $C=0$, then it is called {\it $\lambda$-bilipschitz}. A curve $\gamma:I\to X$ is called a $(\lambda,C)$-{\it  quasigeodesic} if $\gamma$ is a $(\lambda,C)$-quasi-isometric map.

As in \cite{BHK}, a homeomorphism $f:$ $(X,d)\to (X',d')$ is said to be {\it $\theta$-quasisymmetric} if there is a homeomorphism $\theta : [0,\infty) \to [0,\infty)$ such that $d(x,y)\leq t d(x,z)$ implies
$$d'(f(x),f(y)) \leq \theta(t) d'(f(x),f(z))$$
for each $t>0$ and for each triplet $(x, y, z)\in X^3$.

\begin{defn}
Let $\lambda\geq 1$ and $\tau\in (0,1)$. A homeomorphism $f:$ $(X,d)\to (X',d')$ between incomplete and connected metric spaces is a {\it quasisimilarity}, with data $(\theta, \lambda, \tau)$, if

(QS-1) $f$ is $\theta$-quasisymmetric, and

(QS-2) for any $x\in X$, there is a $c_x>0$ (a constant may depend on $x$) such that for all $z,y\in B(x, \tau d(x))$,
$$\lambda^{-1}c_x d(z,y)\leq d'(f(z),f(y))\leq \lambda c_x d(z,y).$$
\end{defn}
Moreover, if there exists a bilipschitz map between two metric spaces, then we say that these two spaces are bilipschitz to each other. If there exists a quasisimilarity between two metric spaces, then we say that these two spaces are quasisimilar to each other.

 \subsection{PQ-isometric maps}
The concepts of cross-difference and strongly PQ-isometric maps were introduced by Buyalo and Schroeder in \cite[Chapter 4]{BuSc}. Let $(X,d)$ be a metric space. For $(x,y,z,u)\in X^4$, the {\it cross-difference} is given by
$$\langle x,y,z,u\rangle=2^{-1}\big(d(x,z)+d(y,u)-d(x,y)-d(z,u)\big).$$

Let $\lambda\geq 1$ and $C\geq 0$. We say that a map $f\colon$ $(X,d)\to (X',d')$ between metric spaces is {\it strongly $\rm{PQ}$-isometric} if for all quadruples $(x,y,z,u)\in X^4$ with $\langle x,y,z,u\rangle\geq 0$,
$$\lambda^{-1}\langle x,y,z,u\rangle-C\leq \big\langle f(x),f(y),f(z),f(u)\big\rangle \leq \lambda\langle x,y,z,u\rangle +C.$$
Moreover, one observes that $\langle x,y,z,u\rangle=-\langle x,z,y,u\rangle$. Therefore, we can write the condition to be a strongly PQ-isometric map as
$$-\theta\big(-\langle x,y,z,u\rangle\big)\leq \big\langle f(x),f(y),f(z),f(u)\big\rangle\leq \theta\big(\langle x,y,z,u\rangle\big),$$
where now  $(x,y,z,u)\in X^4$ is an arbitrary quadruple and $\theta:\mathbb{R}\to \mathbb{R}$ is the control function with
$\theta(t)=(\lambda t)\vee (t/\lambda)+C$.

\section{Gromov hyperbolic spaces}\label{sec-3}

\subsection{Gromov hyperbolicity}\label{sec-g}
Let $(X,d)$ be a metric space and $w\in X$. For $x,y\in X$,  the number
$$(x|y)_w=\langle x,y,w,w\rangle=2^{-1}(d(x,w)+d(y,w)-d(x,y))$$
is called the  {\it Gromov product} of $x,y$ with respect to $w$. One easily finds that
$$\langle x,y,z,u\rangle=-(x|z)_o-(y|u)_o+(x|y)_o+(z|u)_o,$$
for all quadruples $(x,y,z,u)\in X^4$ and for any fixed $o\in X$. Let $\delta\geq 0$.
We say that $X$ is {\it Gromov $\delta$-hyperbolic}, if
$$(x|y)_w\geq (x|z)_w \wedge(z|y)_w-\delta$$ for all $x,y,z,w\in X$. Namely, the triplet $((x|y)_w,(x|z)_w,(y|z)_w)$ is a {\it$\delta$-triple}.

There is an equivalent definition for Gromov hyperbolic geodesic spaces which is known as the Rips condition: Each point on the edge of any geodesic triangle is within the distance $\delta'$ of some point on one of the other two edges, where $\delta'\geq 0$. A geodesic triangle is a set $\Delta=\Delta(x_1, x_2, x_3)=[x_1, x_2]\cup [x_2, x_3]\cup [x_3, x_1]\subset X$. Moreover, $[x_1, x_2]_{\Delta}$ means the edge between $x_1$ and $x_2$ in $\Delta$.
The following result is due to Bonk \cite[Lemma $1.3$]{Bo}. For the definition of tripod maps, see  \cite[p. 284]{Bo}.

\begin{lem}\label{z1}
Let $\Delta=\Delta(x_1,x_2,x_3)$ be a geodesic triangle in a metric space $(X,d)$. If $\Delta$ satisfies the Rips condition with constant $\delta$, then $\Delta$ is $4\delta$-thin. That is, there exists a tripod map $f:\Delta\rightarrow T$ with the following property: If $u,v\in \Delta$ and $f(u)=f(v)$, then $ d(u,v)\leq 4\delta$.
\end{lem}

Suppose that $(X, d)$ is a $\delta$-hyperbolic space. A sequence $\{x_i\}$ in $X$ is called a {\it Gromov sequence} if $(x_i|x_j)_w\rightarrow \infty$ as $i,$ $j\rightarrow \infty.$ Two such sequences $\{x_i\}$ and $\{y_j\}$ are said to be {\it equivalent} if $(x_i|y_i)_w\rightarrow \infty$ as $i\to\infty$.
The {\it Gromov boundary} $\partial_\infty X$ of $X$ is defined to be the set of all equivalence classes of Gromov sequences, and $X^*=X \cup \partial_\infty  X$ is called the {\it Gromov closure} of $X$. See \cite{BS,BrHa,BuSc,Vai-0} for more information.

\begin{lem}\label{b-2} $($\cite[Chapter III.H, Lemmas 3.1 and 3.2]{BrHa}$)$ Suppose X is a proper geodesic $\delta$-hyperbolic space. For each $w\in X$ and $\xi\in \partial_\infty X$, there is a geodesic ray $\gamma:[0,\infty)\to X$ with $\gamma(0)=w$ and $\gamma(\infty)=\xi$. Similarly, for all distinct points $\xi,\eta\in \partial_\infty X$, there is a geodesic line $\gamma:\mathbb{R}\to X$ with $\gamma(-\infty)=\xi$ and $\gamma(\infty)=\eta$.
\end{lem}

For $x\in X$ and $\xi\in \partial_\infty X$, the Gromov product $(x|\xi)_w$ of $x$ and $\xi$ is defined by
$$(x|\xi)_w= \inf \big\{\liminf_{i\rightarrow \infty}(x|y_i)_w\;|\; \{y_i\}\in \xi\big\}.$$
For $\xi,$ $\zeta\in \partial_\infty X$, the Gromov product $(\xi|\zeta)_w$ of $\xi$ and $\zeta$ is defined by
$$(\xi|\zeta)_w= \inf \big\{\liminf_{i\rightarrow \infty}(x_i|y_i)_w\;|\; \{x_i\}\in \xi,\;\{y_i\}\in \zeta\big\}.$$

\begin{lem}\label{z0}$($\cite[Lemma $5.11$]{Vai-0}$)$
Let $X$ be a $\delta$-hyperbolic space with $z,w\in X$, and let $\xi,\xi'\in\partial_\infty X$. Then for any sequences $\{y_i\}\in \xi$, $\{y_i'\}\in \xi'$, we have
\begin{enumerate}
\item  $(z|\xi)_w\leq  \liminf\limits_{i\rightarrow \infty} (z|y_i)_w \leq  \limsup\limits_{i\rightarrow \infty} (z|y_i)_w\leq (z|\xi)_w+\delta;$
\item  $(\xi|\xi')_w\leq  \liminf\limits_{i\rightarrow \infty} (y_i|y_i')_w \leq  \limsup\limits_{i\rightarrow \infty} (y_i|y_i')_w\leq (\xi|\xi')_w+2\delta.$
\end{enumerate}
\end{lem}

\subsection{Busemann functions}
Let $(X,d)$ be a $\delta$-hyperbolic space with $o\in X$ and $\xi\in\partial_\infty X$.
Let $\mathcal{B}(\xi)$ be the class of Busemann functions based at $\xi$ as introduced in~\cite[Section 3.1]{BuSc}. Let $b\in \mathcal{B}(\xi)$ be a Busemann function.
For all $x\in X$, one has
$$b(x)=b_{\xi,o}(x)=b_\xi(x,o)=(\xi|o)_x-(\xi|x)_o.$$
It is known that (cf. \cite[Proposition 3.1.5(1)]{BuSc}),
\be\label{b-1} |b(x)-b(y)|\leq d(x,y)+10\delta, \,\ \mbox{for all } x,y\in X.\ee
According to \cite[Lemma $3.1.1$]{BuSc}, we see that
\be\label{zz0} b_\xi(x,o)\doteq_{2\delta} \{(z_i|o)_x-(z_i|x)_o\}_i =\{d(x,z_i)-d(o,z_i)\}_i\ee
for every Gromov sequence $\{z_i\}\in\xi$.

The Gromov product of $x,y\in X$ with respect to the Busemann function $b=b_{\xi,o}\in \mathcal{B}(\xi)$ is defined by
$$(x|y)_b=2^{-1}(b(x)+b(y)-d(x,y)).$$
By \cite[(3.2) and Example 3.2.1]{BuSc}, we know that
\be\label{zz0.1} (x|y)_b\doteq_{10\delta} (x|y)_o-(x|\xi)_o-(y|\xi)_o.\ee
Similarly, for $x\in X$ and $\zeta\in \partial_\infty X\setminus\{\xi\}$, the Gromov product $(x|\zeta)_b$ of $x$ and $\zeta$ based at $b$ is defined by
$(x|\zeta)_b= \inf \left\{ \liminf_{i\rightarrow \infty}(x|z_i)_b\; |\; \{z_i\}\in \zeta\right\}.$
For a pair of distinct points $\xi_1,\xi_2\in  \partial_\infty X\setminus\{\xi\}$, we define their Gromov product based at $b$ by
$(\xi_1|\xi_2)_b=\inf\left\{\liminf_{i\to\infty} (x_i|y_i)_b\;|\; \{x_i\}\in\xi_1 , \{y_i\}\in\xi_2\right\}.$


\begin{lem}\label{z8}
Let $f:X\to X'$ be a $(\lambda,C)$-quasi-isometry between proper geodesic $\delta$-hyperbolic spaces with $\xi\in \partial_\infty X$ and $f(\xi)=\xi'\in\partial_\infty X'$. Let $b=b_{o,\xi}\in \mathcal{B}(\xi)$ and $b'\in \mathcal{B}(\xi')$. Then there is a control function $\theta:\mathbb{R}\to \mathbb{R}$ with $\theta(t)=(\lambda t)\vee(t/\lambda)+C'$ depending only on $\lambda, C$, and $\delta$ such that, for all $x,y,z,u\in X$,
\be\label{ss1} \langle x',y',z',u'\rangle\leq \theta(\langle x,y,z,u\rangle), \ee
where $f(p)=p'$ for all $p\in X$. Moreover, we have
\be\label{ss2} (x'|z')_{b'}-(x'|y')_{b'}\leq \theta((x|z)_{b}-(x|y)_{b}).\ee
\end{lem}
\bpf The first statement follows from \cite[Theorem $4.4.1$]{BuSc} with $\theta_0(t)=(\lambda t)\vee(t/\lambda)+C_1$ depending only on $\lambda, C$, and $\delta$. It remains to check (\ref{ss2}).

Fix $x,y,z\in X$. For any $\{u_n\}\in \xi$, we have $\{u_n'\}\in\xi'=f(\xi)$ by
\cite[Proposition~6.3]{BS}. Fix a base point $o\in X$. It follows from (\ref{zz0.1}) and Lemma \ref{z0} that
\beq\nonumber
(x|z)_{b}-(x|y)_{b} &\doteq_{20\delta}& (x|z)_o-(x|\xi)_o-(z|\xi)_o-(x|y)_o+(x|\xi)_o+(y|\xi)_o
\\ \nonumber&=& (x|z)_o-(z|\xi)_o-(x|y)_o+(y|\xi)_o
\\ \nonumber&\doteq_{4\delta}&  \{ (x|z)_o-(z|u_n)_o-(x|y)_o+(y|u_n)_o\}_n
\\ \nonumber&=&  \{ \langle x,z,y,u_n\rangle\}_n,
\eeq
and similarly,
$$(x'|z')_{b'}-(x'|y')_{b'} \leq \{ \langle x',z',y',u'_n\rangle\}_n+24\delta.$$
These two estimates, together with (\ref{ss1}), imply that
\beq\nonumber
(x'|z')_{b'}-(x'|y')_{b'} &\leq& \{ \langle x',z',y',u'_n\rangle\}_n+24\delta
\\ \nonumber&\leq& \{\theta_0( \langle x,z,y,u_n\rangle)\}_n+24\delta
\\ \nonumber&\leq&\theta_0((x|z)_{b}-(x|y)_{b}+24\delta)+24\delta
\\ \nonumber&\leq&\theta((x|z)_{b}-(x|y)_{b}),
\eeq
where $\theta(t)=(\lambda t)\vee(t/\lambda)+C'$ with $C'=C_1+24\lambda\delta+24\delta$. Hence  (\ref{ss2}) follows.
\epf

\subsection{Rough starlikeness}\label{sub-rs}
Let $(X,d)$ be a proper geodesic $\delta$-hyperbolic space, and $K\geq 0$. As in \cite{Vai},
we say that $X$ is {\it $K$-roughly starlike} with respect to $\xi$ if for each $x\in X$, there is  a point $\zeta\in\partial_\infty X$ and a geodesic line $\gamma=[\xi,\zeta]$ connecting $\xi$ and $\zeta$ such that  $\dist(x,\gamma)\leq K.$

Notice that a roughly starlike Gromov hyperbolic space contains at least two points on its Gromov boundary. The concept of rough starlikeness with respect to points within the spaces was introduced by Bonk, et al. in \cite{BHK}, which is equivalent to the visual property defined in \cite{BS}. The class of Gromov hyperbolic spaces that are roughly starlike, is very large. For example, it includes metric trees, Gromov hyperbolic domains in $\mathbb{R}^n$ or annular convex metric spaces \cite{HSX,Vai,ZR}, and hyperbolic fillings \cite{BS,BuSc}.
This notation is useful in many studies (cf. \cite{HSX,Vai,ZLR,ZR2}).

In \cite{BHK}, it was shown that bounded uniform spaces are Gromov hyperbolic in the quasihyperbolic metric and roughly starlike with respect to points within the spaces.
In the following, our goal is to show that unbounded uniform spaces are roughly starlike with respect to all points on their Gromov boundaries, see Lemma \ref{z13}.

We begin with some preparations.
Let $(\Omega, d)$ be an incomplete and connected metric space, and $0<\nu\leq 1/2$. Here, $\Omega$ is assumed to be connected in order to avoid the situation that $d(x)=0$ for $x\in \Omega$.
Following the notation of \cite[Chapter 7]{BHK}, a point $x_0$ in $\Omega$ is said to be a {\it $\nu$-annulus point} of $\Omega$, if there is a point $x_1\in\partial\Omega$ such that for $t=d(x_0,x_1)=d(x_0)=\dist(x_0,\partial\Omega),$
the annulus $\{z\in \overline{\Omega}\;|\;\nu t<d(z, x_1)<t/\nu\}$ is contained in $\Omega$. If $x_0$ is not a $\nu$-annulus point of $\Omega$, then it is a {\it $\nu$-arc point} of $\Omega$.
We remark that these concepts are useful in establishing Gromov hyperbolic characterization of uniform domains in annular convex metric spaces, see \cite{HSX}.

\begin{lem}\label{z12}
Let $(\Omega,d)$ be an unbounded $A$-uniform metric space. If $x_0$ is a $\nu$-arc point of  $\Omega$, then there is a quasigeodesic line $\gamma$ in $(\Omega,k)$ with $x_0\in\gamma$.
\end{lem}

\bpf Fix $x_0\in\Omega$ and choose $x_1\in\partial\Omega$ with $d(x_1,x_0)=d(x_0)$. Because  $x_0$ is a $\nu$-arc point, there is a point $x_2\in\partial \Omega$ such that
\be\label{r-2} \nu d(x_0)\leq d(x_1,x_2)\leq \frac{1}{\nu}d(x_0).\ee

As $\Omega$ is $A$-uniform, we observe from \cite[Theorem 3.6]{BHK} that $(\Omega,k)$ is a proper geodesic $\delta$-hyperbolic space with $\delta=\delta(A)\geq 0$.
Choose two quasihyperbolic geodesic rays $\alpha_1$ and $\alpha_2$ joining $x_0$ to $x_1$ and $x_2$, respectively. By \cite[Proposition 3.12]{BHK}, it follows that $\alpha_1$ and $\alpha_2$ are $B$-uniform arcs with $B=B(A)\geq 1$. Next, we proceed to prove the following:

\emph{Claim.}  For each $u\in\alpha_1$, $\ell(\alpha_1[x_1,u])\leq 2B^2d(u)$ and for every $z\in\alpha_2$,
$$\ell(\alpha_2[x_2,z])\leq 2B^2(1+1/\nu)d(z).$$

We only verify the required estimate for the case $z\in \alpha_2$, because the assertion for $u\in \alpha_1$ follows from a similar  argument. To this end, let $y_0$ be the point bisecting the length of $\alpha_2,$ and let $z\in\alpha_2$ be given. We consider the cases corresponding to possible locations of $z$, $y_0$, and $x_0$ in the curve $\alpha_2$.

If $z\in\alpha_2[x_2,y_0]$, then the desired estimate follows from the fact that $\alpha_2$ is $B$-uniform. Thus we are left with the case that $z\in\alpha_2[x_0,y_0]$. We first show that
\be\label{q-2} d(x_0)\leq 2Bd(z). \ee
Indeed, if $d(z,x_0)\leq  d(x_0)/2$, then
$d(z)\geq d(x_0)-d(z,x_0)\geq d(x_0)/2.$
If $d(z,x_0)> d(x_0)/2$, by the uniformity of $\alpha_2$, we obtain that
$$d(z)\geq \frac{1}{B}\ell(\alpha_2[z,x_0])\geq \frac{1}{B}d(z,x_0)\geq \frac{1}{2B}d(x_0).
$$
Thus (\ref{q-2}) is true.

Furthermore, because $\alpha_2$ is $B$-uniform, by (\ref{r-2}) and (\ref{q-2}), it follows that
\beq\label{r-3}
\ell(\alpha_2[x_2,z])&\leq& \ell(\alpha_2[x_2,x_0]) \leq  Bd(x_0,x_2)
\\ \nonumber&\leq&  B(d(x_0,x_1)+d(x_1,x_2))
\\ \nonumber&\leq& B\Big(1+\frac{1}{\nu}\Big)d(x_0)
\leq 2B^2\Big(1+\frac{1}{\nu}\Big)d(z).
\eeq
This proves the claim.

Let $\gamma=\alpha_1 \cup \alpha_2$. It remains to show that $\gamma$ is a quasigeodesic line in $(\Omega,k)$. That is, we need to find some positive numbers $A_1$ and $A_2$ depending only on $A$ and $\nu$ such that
$$k(x,y)\leq \ell_k(\gamma[x,y])\leq A_1 k(x,y)+A_2, \quad \mbox{for all } x,y\in \gamma,$$
where $\ell_k(\gamma[x,y])$ denotes the quasihyperbolic length of the curve $\gamma[x,y]$. Fix $x,y\in \gamma$. Choose two points $y_1\in\alpha_1$ and $y_2\in\alpha_2$ such that
$$
\ell(\alpha_1[x_1,y_1])=\ell(\alpha_2[x_2,y_2])={\nu}d(x_0)/3.
$$
We consider four cases. In what follows, for each $i=1,2,3,4$, $C_i$ is a constant depending only on $A$ and $\nu$.

\begin{enumerate}[(i)]
\item If $x,y\in \alpha_1$ or $x,y\in \alpha_2$, then $\gamma[x,y]$ is evidently a quasihyperbolic geodesic.
\item If $x\in\alpha_1$ and $y\in\alpha_2[y_2,x_0]$, then by (\ref{r-3}), we have
 $$d(x_0,y_2)\leq \ell(\alpha_2)\leq Bd(x_0,x_2)\leq B\Big (1+\frac{1}{\nu}\Big)d(x_0).$$
Moreover, by the claim, we see that
$$d(y_2)\geq \frac{\ell(\alpha_2[x_2,y_2])}{ 2B^2\big (1+1/\nu\big )}=\frac{\nu^2}{6B^2(1+\nu)}d(x_0).$$
Combined with \cite[$(2.16)$]{BHK}, these two estimates ensure that
$$k(x_0,y)\leq k(x_0,y_2)\leq 4A^2\log\left(1+\frac{d(x_0,y_2)}{d(x_0)\wedge d(y_2)}\right)\leq C_1,$$
which guarantees that
$$\ell_k(\gamma[x,y])=  k(x,x_0)+k(x_0,y)\leq 2k(y,x_0)+k(x,y)\leq 2C_1+k(x,y).$$

\item If $x\in\alpha_1[y_1,x_0]$ and $y\in\alpha_2$, then a similar argument as above shows that there is a positive constant $C_2$ such that $\ell_k(\gamma[x,y])\leq 2C_2+k(x,y).$

\item If $x\in\alpha_1[y_1,x_1]$ and $y\in\alpha_2[y_2,x_2]$, then by our choices of the points $y_1$ and $y_2$, we obtain by (\ref{r-2}) that
$$d(x,y)\geq d(x_1,x_2)-d(x_1,x)-d(x_2,y)\geq \nu d(x_0)/3.
$$
This inequality, together with the claim, implies that
$$\ell(\gamma)=\ell(\alpha_1)+\ell(\alpha_2)\leq B d(x_0)+B\Big(1+\frac{1}{\nu}\Big)d(x_0)\leq \frac{3(2\nu+1)B}{\nu^2}d(x,y).$$
Furthermore, we find from the claim that $\gamma[x,y]$ satisfies the condition \cite[(2.14)]{BHK} in \cite[Lemma 2.13]{BHK}. Therefore, it follows from \cite[Lemma 2.13]{BHK} that
\beq\nonumber
\ell_k(\gamma[x,y])&\leq&  8B^2\Big(1+\frac{1}{\nu}\Big)\log \left(1+\frac{\ell(\gamma[x,y])}{d(x)\wedge d(y)}\right)
\\ \nonumber&\leq& 8B^2\Big(1+\frac{1}{\nu}\Big)\log \left(1+\frac{3(1+2\nu)B}{\nu^2}\frac{d(x,y)}{d(x)\wedge d(y)}\right)
\\ \nonumber&\leq& C_3\log\left(1+\frac{d(x,y)}{d(x)\wedge d(y)}\right)+C_4
\\ \nonumber&\leq&  C_3k(x,y)+C_4
\eeq
with $C_3$ and $C_4$ depending only on $B$ and $\nu$.
\end{enumerate}
\vspace{-.5cm}
\epf

\begin{lem}\label{z13}
Let $(\Omega,d)$ be an unbounded $A$-uniform space. Then $(\Omega,k)$ is $K$-roughly starlike with respect to each point of $\partial_\infty(\Omega,k)$ with $K=K(A)$.
\end{lem}

\bpf By \cite[Theorem 3.6]{BHK}, $(\Omega,k)$ is a proper geodesic $\delta$-hyperbolic space with $\delta=\delta(A)\geq 0$. It follows from \cite[Proposition 3.12]{BHK} that there is a natural identification between the one-point extended metric boundary $\partial\Omega\cup\{\infty\}$ and the Gromov boundary $\partial_\infty(\Omega,k)$. According to Lemma \ref{b-2}, for any pair of points on the boundary $\partial \Omega$, there is a quasihyperbolic geodesic line connecting them.

Fix $\omega\in\partial\Omega\cup\{\infty\}$ and $x_0\in\Omega$. We only need to show that there is a quasihyperbolic geodesic line $\gamma$ emanating from $\omega$ such that the quasihyperbolic distance
\be\label{r-4} \dist_k(x_0,\gamma)\leq K\ee
for some constant $K=K(A)\geq 0$. We consider two cases.

Suppose first that $x_0$ is a $\nu$-arc point with $\nu=1/3$.
By Lemma \ref{z12}, there is a $(\lambda, C)$-quasigeodesic $\alpha$ in $(\Omega,k)$ ending at $u,v\in\partial\Omega$ such that $x_0\in\alpha$, where $\lambda$ and $C$ depend only on $A$. Moreover, by the extended stability of Gromov hyperbolic spaces (cf. \cite[Theorem 6.32]{Vai-0}), there is a quasihyperbolic geodesic line $\alpha_0$ joining $u$ and $v$ such that the Hausdorff distance
$$k_{\mathcal{H}}(\alpha,\alpha_0)\leq C_1=C_1(\lambda,C,\delta)=C_1(A),$$
where $C_1(A)\geq 0$. So there is a point $x\in\alpha_0$ such that $k(x_0,x)\leq C_1.$

If $u=\omega$ or $v=\omega$, then (\ref{r-4}) is valid with $K=C_1$. Thus we may assume that $\omega\not\in\{u,v\}$. It follows from \cite[Proposition 2.2]{CDP}  that the extended quasihyperbolic geodesic triangles are $\delta'$-thin for some $\delta'=\delta'(A)\geq 0$. This implies
$$\dist_k(x,\beta_1\cup \beta_2)\leq \delta',$$
where $\beta_1$ and $\beta_2$ are two quasihyperbolic geodesic lines joining $\omega$ to $u$ and $v$, respectively. Therefore, we obtain
$$\dist_k(x_0,\beta_1\cup \beta_2)\leq C_1+\delta':=K.$$
If we choose $\alpha=\beta_1$ or $\alpha=\beta_2$, then the desired estimate (\ref{r-4}) follows in this case.

Next, we assume that $x_0$ is a $\nu$-annulus point with $\nu=1/3$.
By the definition of annulus points, there is a point $x_1\in \partial\Omega$ with $d(x_0)=d(x_0,x_1)$ such that the annulus $\{y\in \overline{\Omega}\;|\;d(x_0)/3<d(y, x_1)<3 d(x_0)\}$ is contained in $\Omega$. Moreover, we claim that
\be\label{t-3} {2d(x_0,x_1)}/{3}\leq d(z)\leq d(x_0,x_1),\ee
for all $z\in S(x_1,d(x_0,x_1))=\{y\in \Omega\; |\;d(y,x_1)=d(x_0,x_1)\}$.

Indeed, for all $u\in\partial \Omega$, we have $d(u,x_1)\geq 3d(x_0,x_1)$ or $d(u,x_1)\leq d(x_0,x_1)/3$. If $d(u,x_1)\geq 3d(x_0,x_1)$, then we get
$$d(z,u)\geq d(u,x_1)-d(z,x_1)\geq 2d(x_0,x_1).
$$
If $d(u,x_1)\leq d(x_0,x_1)/3$, then we obtain
$$d(z,u)\geq d(z,x_1)-d(u,x_1)\geq {2d(x_0,x_1)}/{3},
$$
which yields (\ref{t-3}).

Now, we consider two subcases. If $\omega\in B\big(x_1, d(x_0,x_1)/3\big)$, then there is a quasihyperbolic geodesic line $\gamma$ joining $\omega$ to the infinity point $\infty$. Note that the set $S(x_1,d(x_0, x_1))\cap \gamma$ is non-empty.
To prove (\ref{r-4}), we let $y_0\in S\big(x_1,d(x_0,x_1)\big)\cap \gamma$.
By (\ref{t-3}), we have $d(x_0,y_0)\leq 2d(x_0,x_1)$ and $d(x_0)\wedge d(y_0)\geq {2d(x_0,x_1)}/{3}.$
Thus it follows from \cite[$(2.16)$]{BHK} that
$$k(x_0,y_0)\leq 4A^2\log\left(1+\frac{d(x_0,y_0)}{d(x_0)\wedge d(y_0)}\right)\leq 8A^2=:K
$$
which proves (\ref{r-4}).
If $\omega\in\partial\Omega \cup\{\infty\}\setminus B\big(x_1, d(x_0,x_1)/3\big)$, then a similar argument as above gives that there is a quasihyperbolic geodesic line $\gamma$ joining $\omega$ to $x_1$ with $\dist_k(x_0,\gamma)\leq K$.
\epf

\begin{thm}\label{thm-3}
If $(\Omega,d)$ is an unbounded locally compact uniform metric space, then $(\Omega,k)$ is a proper and geodesic Gromov hyperbolic space which is roughly starlike with respect to any point on the Gromov boundary of $(\Omega,k)$.
\end{thm}
Theorem \ref{thm-3} immediately follows from Lemma \ref{z13} and the first part of \cite[Theorem 3.6]{BHK}. It turns out that, for a given unbounded locally compact uniform space, there is a proper geodesic Gromov hyperbolic space which is roughly starlike with respect to any point on the Gromov boundary.

\section{Uniformization of Gromov hyperbolic spaces}\label{sec-4}
Our goal of this section is to prove Theorem \ref{z3}.  Let $(X,d)$ be a proper geodesic Gromov $\delta$-hyperbolic space with $\delta\geq 0$. Fix $\xi\in\partial_\infty X$ and a Busemann function $b=b_{\xi,o}:X\to \mathbb{R}$ based at $\xi$ with $o\in X$. Consider the family of conformal deformations of $X$ induced by the densities
\be\label{b-0} \rho_\kappa(z)=e^{-\kappa b(z)},\;\;\;\mbox{for}\;\kappa>0.\ee
For $x, y\in X$, we define
$$d_\kappa(x,y)=\inf \int_\gamma \rho_\kappa \,ds,$$
where the infimum is taken over all rectifiable curves $\gamma$ in $(X,d)$ joining the points $x$ and $y$. Thus $d_\kappa$ is a metric on $X$
and we denote the resulting metric spaces by $X_\kappa=(X,d_\kappa)$.
Let $\overline{X_\kappa}$ and
$\partial_\kappa X:=\partial X_\kappa=\overline{X_\kappa}\setminus X_\kappa$
be the metric completion and boundary of $X_\kappa$, respectively. For a rectifiable curve $\gamma$ in $X$, the length of $\gamma$ in the metric $d_\kappa$ is denoted by $\ell_\kappa(\gamma)$.

By (\ref{b-1}), we have
the following {\it Harnack type inequality}: For all $x,y\in X$,
\be\label{q-3} e^{-10\kappa\delta}e^{-\kappa d(x,y)}\leq \frac{\rho_\kappa(x)}{\rho_\kappa(y)}\leq e^{10\kappa\delta}e^{\kappa d(x,y)}.\ee
Applying (\ref{q-3}) and \cite[Theorem 5.1]{BHK}, we obtain the following Gehring-Hayman Theorem. For more background information, we refer to \cite{BHK,RSZ,ZLR} and the references therein.

\begin{thm}\label{q-5} Let $X$ be a geodesic $\delta$-hyperbolic metric space, and let $\rho_\kappa(x)=e^{-\kappa b(x)}$ for all $x\in X$ and $\kappa>0$. There exist $\kappa_0=\kappa_0(\delta)>0$ and $M=20e^{20\kappa\delta}$ such that, if $\kappa\leq \kappa_0$, then
$\ell_\kappa([x,y])=\int_{[x,y]}\rho_\kappa \,ds\leq M d_\kappa(x,y)$
for each geodesic $[x,y]$ in $X$.
\end{thm}

To show Theorem \ref{z3}, we still need an auxiliary result which will be used in the next section.

\begin{lem}\label{z2}
Let $(X,d)$ be a proper geodesic $\delta$-hyperbolic space with $w\in X$, and $\Delta=\Delta(x,y,\xi)$ an extended geodesic triangle with $x,y\in X$ and $\xi\in\partial_\infty X$. If $b=b_{\xi,w}\in \mathcal{B}(\xi)$ and there is a point $w_y\in[x,y]$ satisfying $d(y,w_y)=(x|\xi)_y$, then
\begin{enumerate}
\item\label{q-1}
$b(u)-b(z)\geq d(u,z)-16\delta,\ \mbox{ for any $z\in [x,w_y]$ and for all $u\in [x,z]$,}$
\item\label{z4} $ b(w_y) \doteq_{16\delta} (x|y)_b$.
\end{enumerate}
\end{lem}

\bpf We first prove (\ref{q-1}). For any $z\in [x,w_y]$ and for each $u\in [x,z]$, let $w_x\in[x,y]$ be a point such that $d(x,w_x)=(y|\xi)_x$. According to Lemma \ref{z0}, we see that
$(x|\xi)_y+(y|\xi)_x\leq d(x,y),$
which implies
\be\label{eq-ss0a} w_x\in[x,w_y].\ee
Moreover, again by Lemma \ref{z0}, it follows that
\be\label{ss0} d(w_x,w_y)= d(x,y)-(x|\xi)_y-(y|\xi)_x\leq 2\delta.\ee

Next, consider a sequence of points $a_n\in[x,\xi]$ with $\{a_n\}\in\xi$. Thus for $n$ sufficiently large, by Lemma \ref{z0}, without loss of generality, we assume that
$(y|\xi)_x\leq (y|a_n)_x$ and $(x|\xi)_y\leq (x|a_n)_y.$
Let $\Delta_n=\Delta_n(x,y,a_n)$ be a geodesic triangle with the edges $[x,y]_{\Delta_n}=[x,y]_\Delta$ and $[x,a_n]_{\Delta_n}\subset [x,\xi]_\Delta$. We may pick $w_n\in[x,y]$ with $d(x,w_n)=(y|a_n)_x$ and $d(y,w_n)=(x|a_n)_y$. Observe that
$$d(x,w_x)=(y|\xi)_x\leq d(x,w_n)$$ and
$$d(y,w_y)=(x|\xi)_y\leq d(y,w_n).$$
Thus, we have $w_n\in[w_x,w_y]$.

Furthermore, by (\ref{zz0}), we know that
$$b(u)-b(z)\doteq_{4\delta}\{d(a_n,u)-d(a_n,z)\}_n.$$
Therefore, to demonstrate (\ref{q-1}), it suffices to show that
$$d(a_n,u)-d(a_n,z)\geq d(u,z)-12\delta.$$
We consider three cases.

\emph{Case 1.} Suppose that $z\in[x,w_n]$.  A direct application of Lemma \ref{z1} to the geodesic triangle $\Delta_n$ gives that there are two points $u_0,z_0\in [x,a_n]$ such that
$d(x,u_0)=d(x,u),\;\;d(x,z_0)=d(x,z)$, and $d(u,u_0)\vee d(z,z_0)\leq 4\delta$.
Therefore, we obtain
\beq\nonumber
d(a_n,u)-d(a_n,z)&\geq& d(a_n,u_0)-4\delta-d(a_n,z_0)-4\delta
\\ \nonumber &=& d(u_0,z_0)-8\delta=d(u,z)-8\delta.
\eeq

\emph{Case 2.} Suppose that $z,u\in[w_n,w_y]$. By (\ref{ss0}), we have
$$d(u,z)\leq d(w_x,w_y)\leq 2\delta,$$
and therefore,
$$d(a_n,u)-d(a_n,z)\geq-d(u,z)\geq d(u,z)-4\delta,$$
as desired.

\emph{Case 3.} Suppose that $z\in[w_n,w_y]$ and $u\in[x,w_n]$. We find that
\beq\nonumber
d(a_n,u)-d(a_n,z)&\geq& d(a_n,u)-d(a_n,w_n)-d(w_n,z)
\\ \nonumber &\geq& d(u,w_n)-d(w_n,z)-8\delta\;\;\;\;\;\;\;\;\;\;\;\;\;\;\;\;\;\;\;\;\;\;(\mbox{by Case}\;1)
\\ \nonumber &=& d(u,z)-2d(w_n,z)-8\delta
\\ \nonumber &\geq& d(u,z)-12\delta \;\;\;\;\;\;\;\;(\mbox{since}\; d(w_n,z)\leq d(w_x,w_y)\leq 2\delta),
\eeq
which shows the desired conclusion (\ref{q-1}) of Lemma \ref{z2}.

It remains to prove Lemma \ref{z2} (\ref{z4}). By the statement (\ref{q-1}) of Lemma \ref{z2} and by (\ref{b-1}), it follows that
$$d(x,w_y)-16\delta\leq b(x)-b(w_y)\leq d(x,w_y)+10\delta,$$
and a similar argument gives that
$$d(y,w_y)-16\delta\leq b(y)-b(w_y)\leq d(y,w_y)+10\delta.$$
By the above two inequalities, we obtain
$$b(x)+b(y)-2b(w_y)\doteq_{32\delta} d(x,w_y)+d(y,w_y)=d(x,y).$$
This guarantees that $b(w_y) \doteq_{16\delta} (x|y)_b$, completing the proof.
\epf

\subsection{Proof of Theorem \ref{z3}.}
Suppose that $(X,d)$ is a proper geodesic $\delta$-hyperbolic space and $\partial_\infty X$ contains at least two points. Let $d_\kappa$ be the class of metrics on $X$ induced by the densities \eqref{b-0}. The aim is to show that the conformal deformations $X_\kappa=(X,d_\kappa)$ are unbounded locally compact uniform spaces.

To this end, we first demonstrate that the identity map $(X,d)\to (X,d_\kappa)$ is locally bilipschitz.
Fix $z\in X$. For all $x\in B(z,1)$, we have by (\ref{q-3}) that
$$e^{-\kappa-10\kappa\delta}\rho_\kappa(z)\leq \rho_\kappa(x)\leq
e^{\kappa+10\kappa\delta}\rho_\kappa(z).$$
Let $x$ and $y$ be two distinct points in $B(z,1)$. Choose a geodesic $[x,y]$ in $(X,d)$ joining $x$ to $y$. For all $u\in [x,y]$, we observe from (\ref{q-3}) that
$$\rho_\kappa(u)\leq e^{10\kappa\delta}
e^{\kappa d(x,u)}\rho_\kappa(x)\leq
e^{3\kappa+20\kappa\delta}\rho_\kappa(z).$$
This leads to
$$d_\kappa(x,y)\leq \int_{[x,y]}\rho_\kappa(u)\, |du|\leq e^{3\kappa+20\kappa\delta}\rho_\kappa(z) d(x,y).$$
On the other hand, for every rectifiable curve $\alpha$ in $X$ connecting $x$ and $y$, there exists a subcurve $\alpha_0$ of $\alpha$ starting from $x$ such that $\ell(\alpha_0)=d(x,y)$. Again by (\ref{q-3}), we obtain
\beq\nonumber
d_\kappa(x,y) &=& \inf_\alpha \int_\alpha \rho_\kappa \,ds\geq \inf_\alpha \int_{\alpha_0}\rho_\kappa(u)\, |du|
\\ \nonumber&\geq& \int_{0}^{d(x,y)} e^{-10\kappa\delta}\rho_\kappa(x)  e^{-\kappa t}\ dt
\\ \nonumber&\geq&  e^{-10\kappa\delta-2\kappa}\rho_\kappa(x) d(x,y)
\\ \nonumber&\geq& e^{-20\kappa\delta-3\kappa}\rho_\kappa(z) d(x,y).
\eeq

Hence, we see that the identity map $(X,d)\to (X,d_\kappa)$ is locally bilipschitz.
It follows that $(X,d_\kappa)$ is a locally compact and rectifiably connected metric space.

Next, we check that $(X,d_\kappa)$ is unbounded and incomplete.

\emph{Claim.} Let $\gamma:[0,\infty)\to X$ be a geodesic ray with $\gamma(0)=w\in X$ and $\gamma(\infty)=\zeta$. Take a sequence of points $x_n=\gamma(t_n)$ in $X$ with $t_n\to\infty$ as $n\to\infty$.
\begin{enumerate}
  \item If $\zeta=\xi$, then $\{x_n\}$ tends to infinity in $X_\kappa$;
  \item If $\zeta\neq \xi$, then $\{x_n\}$ is a $d_\kappa$-Cauchy sequence and converges to a point in $\partial_\kappa X$.
\end{enumerate}

We first assume $\zeta=\xi$ and show that $d_\kappa(\gamma(0),\gamma(t_n))\to\infty$ as $n\to\infty.$
Take a geodesic ray $\sigma:[0,\infty)\to X$ with $\sigma(0)=o$ and $\sigma(\infty)=\xi$. According to \cite[Lemma III.3.3, p. 428]{BrHa}, there are constants $T_1,T_2>0$ such that for all $t\geq 0$,
\be\label{ss3} d(\gamma(T_1+t),\sigma(T_2+t))\leq 5\delta.\ee
Let $t\geq 2T_1+2T_2$. For $n\geq 1$ sufficiently large, by Lemma \ref{z0}, without loss of generality, we may assume that
\begin{itemize}
\item[(i)] $(\xi|o)_{\gamma(t)}\leq (\sigma(t+T_2-T_1+n)|o)_{\gamma(t)}+2\delta$,
and
\item[(ii)] $(\xi|\gamma(t))_o\geq (\sigma(t+T_2-T_1+n)|\gamma(t))_o-2\delta.$
\end{itemize}
Then we obtain from the above two estimates and (\ref{ss3}) that
\beq\nonumber
b(\gamma(t)) &=& (\xi|o)_{\gamma(t)}-(\xi|\gamma(t))_o
\\ \nonumber&\leq& (\sigma(t+T_2-T_1+n)|o)_{\gamma(t)}-(\sigma(t+T_2-T_1+n)|\gamma(t))_o+4\delta
\\ \nonumber&=& d(\sigma(t+T_2-T_1+n),\gamma(t))-d(o,\sigma(t+T_2-T_1+n))+4\delta
\\ \nonumber&\leq& d(\sigma(t+T_2-T_1+n),\sigma(t+T_2-T_1))+d(\sigma(t+T_2-T_1),\gamma(t))
\\ \nonumber& & ~~~ -(t+T_2-T_1+n)+4\delta
\\ \nonumber&\leq& 9\delta-t-T_2+T_1\leq 9\delta-T_2+T_1.
\eeq
This, together with Theorem \ref{q-5}, shows that
\beq\nonumber
d_\kappa(\gamma(0),\gamma(t_n))&\geq& \frac{1}{20e^{20e^{\kappa\delta}}}\ell_\kappa(\gamma[0,t_n])
\geq \frac{t_n}{20e^{29\kappa\delta+\kappa(T_2-T_1)}}\to\infty,\;\;\mbox{as}\;n\to \infty.
\eeq
Therefore, the first part of the claim is true.

For the other case $\{x_n\}\in\zeta\neq \xi$, we check that $\{x_n\}=\{\gamma(t_n)\}$ is a $d_\kappa$-Cauchy sequence. Because $\partial_\infty X$ contains at least two points, the existence of $\zeta$ is guaranteed.
Again by Lemma \ref{z0}, for $t$ sufficiently large, we obtain that
\be\label{t-4} b(\gamma(t))\geq d(o,\gamma(t))-2(\xi|\gamma(t))_o-2\delta\geq t-d(o,w)-2(\xi|\zeta)_o-8\delta.       \ee
Without loss of generality, we may assume that $t_m\geq t_n>t$. By (\ref{t-4}), we have
\beq\label{t-5}
d_\kappa(\gamma(t_n),\gamma(t_m)) &\leq& \int_{t_n}^{t_m} e^{-\kappa b(\gamma(t))}\, dt
\\ \nonumber&\leq&  e^{\kappa[d(o,w)+2(\xi|\zeta)_o+8\delta]} \int_{t_n}^{t_m} e^{-\kappa t} \,dt
\\ \nonumber&\leq& {\kappa}^{-1} e^{\kappa[d(o,w)+2(\xi|\zeta)_o+8\delta]} e^{-\kappa t_n}.
\eeq
Because $\xi\neq \zeta$, we have $(\xi|\zeta)_o<\infty$. Thus we see from (\ref{t-5}) that $\{x_n\}=\{\gamma(t_n)\}$ is $d_\kappa$-Cauchy. As $\{x_n\}$ is a Gromov sequence in $X$, $\{x_n\}$ is $d_{\kappa}$-convergent to some point in $\partial_{\kappa} X$. Hence, the claim holds.

Set $d_\kappa(x)=\dist_\kappa(x,\partial_\kappa X)$ for all $x\in X$. We  prove that
\be\label{q-4} d_\kappa(x)\geq \frac{1}{2\kappa e^{10\kappa\delta}}\rho_\kappa(x). \ee
For all points $y\in X$ satisfying $d(x,y) \geq 1/\kappa$ and for any rectifiable curve $\gamma$ joining $x$ to $y$, there is a subcurve $\gamma_0$ of $\gamma$ with $x\in\gamma_0$ and $\ell(\gamma_0)=1/\kappa$. Using (\ref{q-3}), we get
\beq\nonumber
d_\kappa(x,y) &=& \inf_\gamma \int_\gamma \rho_\kappa \,ds\geq \inf_\gamma \int_{\gamma_0}\rho_\kappa(z)\, |dz|
\\ \nonumber&\geq& \int_{0}^{1/\kappa} e^{-10\kappa\delta}\rho_\kappa(x)  e^{-\kappa t}\ dt
\\ \nonumber&=& \frac{1-e^{-1}}{\kappa e^{10\kappa\delta}}\rho_\kappa(x),
\eeq
which yields (\ref{q-4}).

Now, we are ready to show that $X_\kappa$ is uniform. For this, it suffices to verify that every geodesic $[x,y]$ in $(X,d)$ is a uniform arc in $X_\kappa$. The quasiconvexity of $[x,y]$ in $X_\kappa$ follows from Theorem \ref{q-5}. It remains to check the second condition of uniformity.

Let $z\in [x,y]$. Consider the extended geodesic triangle $\Delta=[x,\xi]\cup [\xi,y]\cup [y,x].$
Take a point $w_y\in[x,y]$ with $d(y,w_y)=(x|\xi)_y$. We only need to consider the case $z\in [x, w_y]$, because the case $z\in[w_y,y]$ follows from a similar argument by \eqref{eq-ss0a}.

Finally, by Lemma \ref{z2}(\ref{q-1}), we find that
$$b(u)\geq b(z)+d(u,z)-16\delta, \quad \mbox{for each } u\in [x,z].$$
By (\ref{q-4}), we obtain
\beq\nonumber
\ell_\kappa([x,z]) &=&  \int_{[x,z]} \rho_\kappa(u)\, |du|
\\ \nonumber&\leq& \int_{[x,z]}e^{-\kappa b(z)-\kappa d(u,z)+16\delta \kappa}\,|du|
\\ \nonumber&\leq& \rho_\kappa(z)e^{16\delta \kappa}\int_0^\infty e^{-\kappa t}\,dt
= e^{16\delta \kappa}\kappa^{-1}\rho_\kappa(z)
\\ \nonumber&\leq&  2e^{26\kappa \delta}d_\kappa(z).
\eeq

\qed

\section{Uniformization and quasihyperbolization}\label{sec-5}

This section is devoted to the proof of Theorem \ref{main thm-1}.
By Theorems \ref{z3} and \ref{thm-3}, it follows that there is a correspondence between unbounded locally compact uniform spaces and proper geodesic Gromov hyperbolic spaces that are roughly starlike with respect to points on the Gromov boundaries. It remains to show that there is a correspondence of maps between Gromov hyperbolic spaces and maps between uniform spaces.

\textbf{Outline of the proof.} Let $X$ be a proper geodesic Gromov hyperbolic space which is roughly starlike with respect to a point on the Gromov boundary, and let $[X]$ be the set of all proper geodesic Gromov hyperbolic spaces which are bilipschitz to $X$ and roughly starlike with respect to points on their Gromov boundaries. Similarly, let $\Omega$ be an unbounded locally compact  uniform space and $[\Omega]$ the set of all locally compact unbounded uniform spaces which are quasisimilar to $\Omega$.

As in \cite[Chapter 4]{BHK}, let us denote the uniformization procedure $X\to X_\kappa$ by $\mathcal{D}$ for $0<\kappa \leq \kappa_0(\delta)$, where $\mathcal{D}$ stands for {\it dampening}. Similarly, we may {\it quasihyperbolize} an uniform space $\Omega$ by considering its quasihyperbolic metric $k$. Denote this association $\Omega\to (\Omega,k)$ by $\mathcal{Q}$.

Let $\mathcal{G}$ be the bilipschitz classes of proper geodesic roughly starlike Gromov hyperbolic spaces,
and $\mathcal{U}$ the quasisimilarity classes of unbounded locally compact uniform spaces.
To prove Theorem \ref{main thm-1}, we only need to show the following assertions:
\begin{enumerate}[(S-1)]
\item\label{qq-1} $\mathcal{Q}([\Omega])=[(\Omega,k)]\in \mathcal{G}$,
\item\label{qq-2} $\mathcal{D}([X])=[X_\kappa]\in \mathcal{U}$,
\item\label{qq-3} $\mathcal{Q\circ D}([X])=\mathcal{Q} ([X_\kappa])=[(X_\kappa,k_\kappa)]=[X]$,
\item\label{qq-4} $\mathcal{D\circ Q}([\Omega])=\mathcal{D} ([(\Omega,k)])=[\Omega_\kappa]=[\Omega]$,
\end{enumerate}
where $k_\kappa$ is the quasihyperbolic metric of a uniform space $X_\kappa$ and $\Omega_\kappa$ is the conformally deformed space of $(\Omega,k)$ which was introduced in Section \ref{sec-4}. The proof is similar in spirit to that in \cite{BHK}, while our arguments concerning Busemann functions are more complicated and used in a different manner.

\subsection{Auxiliary lemmas} We first establish certain auxiliary results. In this subsection, we assume that $(X,d)$ is a proper and geodesic space that is $\delta$-hyperbolic. Let $X_\kappa=(X,d_\kappa)$ be its uniformization induced by $\rho_{\kappa}(x)=e^{-\kappa b(x)}$ for $0<\kappa\leq \kappa_0(\delta)$, where $b\in \mathcal{B}(\xi)$ and $\xi\in \partial_{\infty} X$.

\begin{lem}\label{z5}
There is a constant $C=C(\delta)\geq 1$ such that
\be\label{zz2} \kappa^{-1} e^{-\kappa (x|y)_b}\big(1\wedge [\kappa d(x,y)]\big)\asymp_C d_\kappa (x,y) ~ \mbox{ for all $x,y\in X$}. \ee
\end{lem}

\bpf
Let $\Delta=\Delta(x,y,\xi)$ be an extended geodesic triangle, and let $w_y\in[x,y]$ be a point such that $d(y,w_y)=(x|\xi)_y$. By Lemma \ref{z2}(\ref{z4}), we see that
$$b(w_y) \doteq_{16\delta} (x|y)_b.$$
To prove (\ref{zz2}),  we only need to check that
\be\label{u-2} d_\kappa (x,y)\asymp_{C(\delta)} \kappa^{-1}{\rho_\kappa (w_y)}\big(1\wedge [\kappa d(x,y)]\big).\ee
Assume first that $\kappa d(x,y)\leq 1$. For all $u\in [x,y]$,
$$\kappa d(u,w_y)\leq \kappa d(x,y)\leq 1.$$
Setting $C_0=e^{10\kappa\delta}$. By (\ref{q-3}), we have
\be\label{u-2b}\frac{1}{C_0e}\rho_\kappa(w_y)\leq \frac{1}{C_0}e^{-\kappa d(w_y,u)}\rho_\kappa(w_y)\leq \rho_\kappa(u)\leq C_0e^{\kappa d(w_y,u)}\rho_\kappa(w_y)\leq C_0e\rho_\kappa(w_y).\ee
This implies that
$$d_\kappa (x,y) \leq \int_{[x,y]}\rho_\kappa(u)\,|du| \leq  C_0e\rho_\kappa(w_y) d(x,y).$$
For the other direction, \eqref{u-2b} together with Theorem \ref{q-5}, guarantees that
$$d_\kappa (x,y) \geq \frac{1}{20C_0^2}\int_{[x,y]}\rho_\kappa(u)\,|du| \geq \frac{\rho_\kappa(w_y)}{20eC_0^3} d(x,y).   $$

Next, we consider the case that $\kappa d(x,y)>1$. By Lemma \ref{z2}(\ref{q-1}) and by symmetry, we have
$$e^{-\kappa b(u)}\leq e^{-\kappa b(w_y)} e^{-\kappa d(u,w_y)}e^{ 16\delta\kappa}$$
for all $u\in [x,y]$.
Hence we obtain
$$\rho_\kappa(u)\leq e^{ 16\delta\kappa} \rho_\kappa(w_y)e^{-\kappa d(u,w_y)},$$
which ensures that
\beq\nonumber
d_\kappa (x,y) &\leq& \int_{[x,y]}\rho_\kappa(u)\,|du|
\\ \nonumber&\leq& e^{ 16\delta\kappa}\rho_\kappa(w_y)\int_{[x,y]}e^{-\kappa d(u,w_y)} \,|du|
\\ \nonumber&\leq& 2e^{ 16\delta\kappa}\rho_\kappa(w_y)\int_0^\infty e^{-\kappa t}\,dt
\\ \nonumber&=& 2e^{ 16\delta\kappa}\kappa^{-1}\rho_\kappa(w_y).
\eeq

Now by (\ref{b-1}), we observe that
$b(u)\leq b(w_y)+d(u,w_y)+10\delta$ for all $u\in[x,y]$.
Therefore,
$\rho_\kappa(u)\geq \rho_\kappa(w_y)e^{-\kappa d(u,w_y)}e^{-10\delta\kappa}.$
By Theorem \ref{q-5}, we obtain
\beq\nonumber
d_\kappa (x,y) &\geq& \frac{1}{20C_0^2} \int_{[x,y]}\rho_\kappa(u)\,|du|
\\ \nonumber&\geq& \frac{\rho_\kappa(w_y)}{20C_0^3}\int_0^{\frac{1}{2}d(x,y)}e^{-\kappa t}\, dt
\\ \nonumber&=& \frac{\rho_\kappa(w_y)}{20C_0^3} \frac{1}{\kappa}(1-e^{-\frac{\kappa}{2}d(x,y)})
\\ \nonumber&>& \frac{\rho_\kappa(w_y)}{20C_0^3} \frac{1}{\kappa}(1-e^{-\frac{1}{2}}),
\eeq
which proves (\ref{u-2}), completing the proof.
\epf

\begin{lem}\label{z6}
There is a natural identification $\phi:\partial_\infty X\to \partial_\kappa X\cup\{\infty\}$.
\end{lem}
\bpf
Firstly, we show that there is a well-defined map $\phi:\partial_\infty X\to \partial_\kappa X\cup \{\infty\}$ with $\phi(\xi)=\infty.$
Let $w\in X$. For any sequence $\{x_n\}\in \xi$, without loss of generality, we may assume that $\kappa d(x_n,w)\geq 1$. Thus by (\ref{zz0.1}), we have
$$(x_n|w)_b\doteq (x_n|w)_o-(x_n|\xi)_o-(\xi|w)_o\to -\infty \mbox{ as } n\to \infty.$$
This, together with (\ref{zz2}), shows that
$$d_\kappa(x_n,w) \asymp \kappa^{-1} e^{-\kappa (x_n|w)_b}\to \infty,$$
as desired. On the other hand, for any sequence $\{x_n\}\in \zeta\neq\xi$, we may check that $\{x_n\}$ is $d_{\kappa}$-convergent to some point in $\partial_\kappa X$. By (\ref{zz0.1}), we get
$$(x_n|x_m)_b\doteq (x_n|x_m)_o-(x_n|\xi)_o-(x_m|\xi)_o\to +\infty \mbox{ as } n,m\to \infty.$$
This, together with (\ref{zz2}), gives that
$$d_\kappa(x_n,x_m)\leq C(\delta)\kappa^{-1} e^{-\kappa (x_n|x_m)_b}\to 0 \mbox{ as } n,m\to \infty.$$
Thus we see that $\{x_n\}$ is a $d_{\kappa}$-Cauchy sequence. Furthermore, from the argument in the proof of Theorem \ref{z3},
we know that the identity map $(X,d)\to (X,d_\kappa)$ is locally bilipschitz. This guarantees that the sequence $\{x_n\}$ is $d_{\kappa}$-convergent to a point in $\partial_\kappa X$.

Next, we check that $\phi$ is injective. Assume that $\phi(\zeta)=\phi(\eta)$ for some $\zeta,\eta\in \partial_\infty X$. If $\phi(\zeta)=\infty$, then for any Gromov sequence $\{x_n\}\in \zeta$, we have $d_\kappa(x_n,w)\to \infty$ as $n\to \infty$. By (\ref{zz0.1}) and (\ref{zz2}), we have
$$(x_n|w)_o-(x_n|\xi)_o-(\xi|w)_o\doteq (x_n|w)_b\to-\infty \ \mbox{\ \mbox{ as $n \to \infty$}.}$$
This implies $(x_n|\xi)_o\to +\infty$ as $n\to \infty$. Therefore, we obtain $\zeta=\xi=\eta$.

Now we assume $\phi(\zeta)=\phi(\eta)\neq \infty$ and thus $\zeta\neq\xi\neq\eta$. For any Gromov sequences $\{x_n\}\in\zeta$ and $\{y_n\}\in\eta$,  we have
$$
	d_\kappa(x_n,y_n)\to 0 \mbox{ as } n\to\infty,
$$
because $\phi(\zeta)=\phi(\eta)\neq\infty$. If $\kappa d(x_n,y_n)\leq 1$, we have
$$
(x_n|y_n)_o\geq (x_n|x_n)_o-d(x_n,y_n)\to+\infty \mbox{ as } n\to\infty,
$$
which gives $\zeta=\eta$.
If $\kappa d(x_n,y_n)> 1$, then by (\ref{zz2}) and (\ref{zz0.1}), we obtain
$$(x_n|y_n)_o-(x_n|\xi)_o-(\xi|y_n)_o\doteq (x_n|y_n)_b\to +\infty\;\;\;\mbox{as}\;n\to\infty.$$
This guarantees  $(x_n|y_n)_o\to \infty$, and thus, $\zeta=\eta$. Therefore, $\phi$ is injective.

It suffices to show that $\phi$ is surjective. For each $a\in \partial_\kappa X$ and for any $d_\kappa$-Cauchy sequence $\{x_n\}$ with $d_\kappa(x_n,a)\to 0$ as $n\to \infty$, we only need to prove that $\{x_n\}$ is a Gromov sequence in the space $(X,d)$, because $\phi(\xi)=\infty$.

By (\ref{zz2}), it follows that there is a constant $C=C(\delta)$ such that
$$d_\kappa(x_n,x_m)\asymp_C  \kappa^{-1} e^{-\kappa (x_n|x_m)_b}\big(1\wedge [\kappa d(x_n,x_m)]\big).$$
Because $d_\kappa(x_n,x_m)\to 0$ as $n,m\to \infty$, there are two possibilities: either $d(x_n,x_m)\to 0$ or $(x_n|x_m)_b\to +\infty$ as $n,m\to \infty$.
If $d(x_n,x_m)\to 0$ as $n,m\to \infty$, then $\{x_n\}$ is $d$-convergent to some point in $X$.
As the identity map $(X,d)\to (X,d_\kappa)$ is locally bilipschitz, we find that $\{x_n\}$ is also $d_{\kappa}$-convergent to a point in $X$, which is a contradiction. Therefore, we have $(x_n|x_m)_b\to +\infty$ as $n,m\to \infty$. By (\ref{zz0.1}), we obtain
$$(x_n|x_m)_o-(x_n|\xi)_o-(x_m|\xi)_o \doteq_{C(\delta)} (x_n|x_m)_b\to +\infty \mbox{ as } n,m\to \infty,$$ which leads to  $(x_n|x_m)_o\to \infty$, completing the proof.
\epf

\begin{lem}\label{z7} If $(X,d)$ is $K$-roughly starlike with respect to $\xi$, then for all $x\in X$,
\be\label{zz3} \frac{1}{2\kappa e^{10\kappa\delta}}\rho_\kappa(x)\leq d_\kappa(x)\leq \frac{e^{16\kappa\delta}\rho_\kappa(x)}{\kappa}(2e^{\kappa K}-1).\ee
\end{lem}

\bpf As the left inequality in $(\ref{zz3})$ follows from (\ref{q-4}), it suffices to verify the right inequality.
Fix $x\in X$. Because $(X,d)$ is $K$-roughly starlike with respect to $\xi$, there is a geodesic line $\gamma=[\xi,\eta]$ such that $\gamma(-\infty)=\xi$, $\gamma(\infty)=\eta$, and
$d(x,y)\leq K$ for some point $y\in[\xi,\eta]$. Moreover, for all $z\in\gamma[y,\eta]$, by \cite[Lemma 3.1.2]{BuSc} and (\ref{zz0}), we see that
\beq\nonumber b(z)-b(x)&\geq & b_\xi(z,y)-b_\xi(x,y)-12\delta
\\ \nonumber&\geq& \{d(z,\gamma(-n))-d(y,\gamma(-n))\}_n-\\ \nonumber
& &\{d(x,\gamma(-n))-d(y,\gamma(-n))\}_n-16\delta\\ \nonumber
&\geq& d(y,z)-d(x,y)-16\delta
\\ \nonumber&\geq& d(y,z)-K-16\delta.
\eeq
This guarantees that
\be\label{h-1}\rho_\kappa(z)\leq e^{16\kappa\delta+\kappa K}\rho_\kappa(x)e^{-\kappa d(y,z)}.\ee

By Lemma \ref{z6}, we know that there is a natural identification between  $\partial_\infty X$ and $\partial_\kappa X\cup\{\infty\}$. Thus, we may regard $\eta$ be a point in $\partial_\kappa X$.
Using (\ref{h-1}) and (\ref{q-3}), we obtain
\beq\nonumber
d_\kappa(x) &\leq& d_\kappa(x,y)+d_\kappa(y,\eta)
\\ \nonumber&\leq& \int_{[x,y]} \rho_\kappa(u)\,|du|+\int_{[y,\eta]} \rho_\kappa(z)\,|dz|
\\ \nonumber&\leq&  e^{16\kappa\delta}\rho_\kappa(x)\left(\int_0^K e^{\kappa t}\,dt+ e^{\kappa K}\int_0^\infty e^{-\kappa t}\,dt\right)
\\ \nonumber&=&  e^{16\kappa\delta}\kappa^{-1}  \rho_\kappa(x) (2e^{\kappa K}-1).
\eeq
\vspace{-.6cm}
\epf

\subsection{Proof of (S-\ref{qq-1})} See \cite[Proposition $4.36$]{BHK}.\qed

\subsection{Proof of (S-\ref{qq-2})}
Let $X$ and $X'$ be  proper geodesic $\delta$-hyperbolic spaces which are $K$-roughly starlike with respect to $\xi\in\partial_\infty X$ and $\xi'\in\partial_\infty X'$, respectively. Denote by $\rho_\kappa(x)=e^{-\kappa b(x)}$ and $\rho_{\kappa'}(z')=e^{-\kappa' b'(z')}$ for all $x\in X$ and $z'\in X'$, where $b\in \mathcal{B}(\xi)$ and $b'\in \mathcal{B}(\xi')$ are  Busemann functions. Let $X_\kappa$ and $X'_{\kappa'}$ be the deformed spaces induced by $\rho_\kappa$ and $\rho_{\kappa'}$, respectively, where $\kappa$ and $\kappa'$ depend only on $\delta$.  Thus, to prove (S-\ref{qq-2}), we only need to show the following result.

\begin{lem}\label{z-1}
Suppose that $f:(X,d)\to (X',d')$ is $M$-bilipschitz with $f(\xi)=\xi'$. Then $f:X_\kappa \to X'_{\kappa'}$ is quasisimilar with the data depending only on $K,$ $M$, and $\delta$.
\end{lem}

\bpf We first check the second requirement (QS-2) for the quasisimilarity of $f$. It follows from Lemma \ref{z5} that there is a constant $C_1=C_1(\delta)$ such that
\be\label{ss4} 1\wedge[\kappa d(x,y)]\asymp_{C_1} \kappa e^{\kappa (x|y)_b}d_\kappa(x,y) \ \mbox{ for all $x,y\in X$}.\ee
By (\ref{b-1}), we see that
\be\label{ss5} (x|y)_b\leq b(x)+5\delta \ \ \mbox{for all $x,y\in X$}. \ee
By Lemma \ref{z7}, it follows that there is a constant $C_2=C_2(\delta, K)$ such that
\be\label{ss6} \kappa d_\kappa (z) \asymp_{C_2} \rho_\kappa(z)=e^{-\kappa  b(z)} \ \ \mbox{ for all $z\in X$}.\ee

Choose $\tau\in(0,1)$ such that
$$\tau'=C_1C_2e^{5\delta \kappa}\frac{2\tau}{1-\tau}\leq 1\wedge \frac{\kappa}{M\kappa'}.$$
For all $z\in X$ and for each pair of points $x, y$ in the metric ball $B_\kappa(z,\tau d_\kappa (z))$ of $(X,d_\kappa)$, we claim that
\be\label{ss7} \kappa d(x,y)\leq 1 \,\;\;\;\mbox{and}\;\;\;\; \kappa' d'(x',y')\leq 1. \ee
We use the notation $f(p)=p'$ for all $p\in X$. Because $x,y\in B_\kappa(z,\tau d_\kappa (z))$,
 $$d_\kappa(x,y)\leq \frac{2\tau}{1-\tau} d_\kappa(x).$$
This, together with (\ref{ss4}), (\ref{ss5}), and (\ref{ss6}), shows that
$$1\wedge[\kappa d(x,y)]\leq C_1\kappa e^{\kappa b(x)+5\delta \kappa}d_\kappa(x,y)\leq C_1C_2 e^{5\kappa \delta} \frac{d_\kappa(x,y)}{d_\kappa(x)}\leq \tau'\leq 1,$$
and, therefore,
$\kappa d(x,y)\leq 1$,
which is the first part of (\ref{ss7}). We also see that
\be\label{ss8} \kappa d(x,y)\leq \tau'.\ee
 Because $f$ is $M$-bilipschitz, we have
$$\kappa' d'(x',y')\leq \kappa'M d(x,y)\leq \kappa'M {\tau'}/{\kappa}\leq 1,$$
by our choice of $\tau'$. This yields (\ref{ss7}).

Moreover, we see from (\ref{b-1}) and (\ref{ss8}) that
\be\label{ss9}\;\; (x|y)_b\doteq_{\tau'/\kappa} \frac{1}{2}[b(x)+b(y)]\doteq_{\tau'/\kappa+10\delta} b(z),\ee
and similarly\;$(x'|y')_{b'}\doteq b'(z')$.
Furthermore, we get
\beq\label{zz4}\;\;\;\;\;\;\;\;\;\;\;
d_\kappa(x,y)&\asymp_{C_1}& {\kappa}^{-1}e^{-\kappa (x|y)_b} \big(1\wedge [\kappa d(x,y)]\big)\;\;\;\;\;\;\;\;\;\;\;\;\;\;(\mbox{by}\;(\ref{ss4}))
\\ \nonumber&\asymp_{e^{2\tau'+10\delta\kappa}}&  e^{-\kappa b(z)} d(x,y)\;\;\;\;\;\;\;\;\;\;\;\;\;\;\;\;\;\;\;\;\;\;\;\;\;\;(\mbox{by}\;(\ref{ss7})\;\mbox{and}\;(\ref{ss9}))
\\ \nonumber&=&\rho_\kappa(z) d(x,y)
\\ \nonumber& \asymp_{C_2} &
\kappa d_\kappa(z) d(x,y).\;\;\;\;\;\;\;\;\;\;\;\;\;\;\;\;\;\;\;\;\;(\mbox{by}\;(\ref{ss6}))
\eeq
Also, a similar argument as above gives that
\be\label{qq-6} d'_{\kappa'}(x',y')\asymp \kappa' d'_{\kappa'}(z') d'(x',y').\ee
As $f:(X,d)\to(X',d')$ is $M$-bilipschitz,  we see from (\ref{zz4}) and (\ref{qq-6}) that the second requirement (QS-2) of quasisimilarity holds.

It remains to show that the induced map $f:X_\kappa \to X'_{\kappa'}$ is quasisymmetric. By \cite[Theorem $6.6$]{Vai-4}, we only need to find some constant $H\geq 1$ depending only on $M,$ $\delta$, and $K$ such that for each triple of distinct points $x,y,z\in X$,
\be\label{s-5} d_\kappa(x,y) \leq d_\kappa(x,z)\;\;\;\;\;\mbox{implies}\;\;\;\;\; d'_{\kappa'}(x',y')\leq H d'_{\kappa'}(x',z'),\ee
because $X_\kappa$ and $X_{\kappa'}$  are both quasiconvex metric spaces.

As $d_\kappa(x,y)\leq d_\kappa(x,z)$, by (\ref{ss4}), we have
$$
e^{\kappa (x|z)_b-\kappa (x|y)_b}\frac{1\wedge [\kappa d(x,y)]}{1\wedge [\kappa d(x,z)]}\leq C_1^2.
$$
Let
$C_3=C_1C_2e^{2\tau'+10\delta\kappa}.$
We divide the proof of (\ref{s-5}) into three cases.

\emph{Case $\textup{I.}$} Suppose that $C_3\kappa d(x,z)<\tau.$
By a similar argument as in (\ref{zz4}), we find that
$$\frac{d_\kappa(x,y)}{d_\kappa(x)}\leq \frac{d_\kappa(x,z)}{d_\kappa(x)}\leq C_3\kappa d(x,z)<\tau,$$
and thus, $y,z\in B_\kappa(x,\tau d_\kappa (x))$. Since $f$ satisfies the second requirement (QS-2) of quasisimilarity, we know that (\ref{s-5}) is true.

\emph{Case $\textup{II.}$} Suppose that  $C_3\kappa d(x,z)\geq\tau$ and $\kappa d(x,y)<1$.
Because $f:(X,d)\to(X',d')$ is $M$-bilipschitz, we observe that
\be\label{s-7} \kappa'd'(x',z')\geq \frac{\tau \kappa'}{C_3 M\kappa}\;\;\;\mbox{and}\;\;\; d'(x',y')\leq \frac{M}{\kappa}.\ee
Now by (\ref{b-1}), we see that
\be\label{s-8}(x'|z')_{b'}-(x'|y')_{b'}=\frac{1}{2}\Big(b'(z')-b'(y')-d'(x',z')+d'(x',y')\Big)\leq d'(x',y')+5\delta.\ee
It follows from Lemma \ref{z5}, (\ref{s-7}), and (\ref{s-8}) that
$$\frac{d'_{\kappa'}(x',y')}{d'_{\kappa'}(x',z')}\leq C_1^2 e^{\kappa'(x'|z')_{b'}-\kappa'(x'|y')_{b'}}\frac{1\wedge [\kappa' d'(x',y')]}{1\wedge [\kappa' d'(x',z')]}\leq H_1$$
for some constant $H_1=H_1(\delta, K, M)$.

\emph{Case $\textup{III.}$} Suppose that $C_3\kappa d(x,z)\geq\tau$ and $\kappa d(x,y)\geq 1$.
Again we have as in (\ref{s-7}),
\be\label{s-9}\kappa'd'(x',z')\geq \frac{\tau \kappa'}{C_3 M\kappa}.\ee
Moreover, by (\ref{ss4}), we see that
$$e^{\kappa(x|z)_{b}-\kappa(x|y)_{b}}\leq C_1^2\frac{d_\kappa(x,y)}{d_\kappa(x,z)}\frac{1\wedge [\kappa d(x,z)]}{1\wedge [\kappa d(x,y)] }\leq C_1^2.$$
This, together with Lemma \ref{z8}, shows that there is a constant $C_4=C_4(C_1,\delta,M)$ such that
\be\label{s-10}(x'|z')_{b'}-(x'|y')_{b'}\leq C_4.\ee
Therefore, by Lemma \ref{z5}, (\ref{s-9}), and (\ref{s-10}), we obtain
$$\frac{d'_{\kappa'}(x',y')}{d'_{\kappa'}(x',z')}\leq C_1^2 e^{\kappa'(x'|z')_{b'}-\kappa'(x'|y')_{b'}} \frac{1\wedge [\kappa' d'(x',y')]}{1\wedge [\kappa' d'(x',z')]}\leq H_2$$
for some constant $H_2=H_2(\delta, K, M)$. This gives (\ref{s-5}), completing the proof.
\epf

\br
Note that the assumption $f(\xi)=\xi'$ is essential. As $f:(X,d)\to(X',d')$ is bilipschitz, it is not difficult to see that $(X',d')$ is $K'$-roughly starlike with respect to $f(\xi)\in\partial_\infty X'$. Moreover, one observes from Lemma \ref{z6} that both $X_\kappa$ and $X'_{\kappa'}$ are unbounded, $\xi$ and $\xi'$ correspond to the infinity point, respectively. Because quasisymmetric maps fix the infinity point between unbounded metric spaces (cf. \cite[Corollary 2.6]{TV}), we see that $f(\xi)=\xi'$, as required.
\er

\subsection{Proof of (S-\ref{qq-3})}
Let $(X,d)$ be a proper geodesic $\delta$-hyperbolic space which is $K$-roughly starlike with respect to $\xi\in\partial_\infty X$, and let $b\in \mathcal{B}(\xi)$ be a Busemann function based at $\xi$, and let $X_\kappa$ be the resulting space of $X$ induced by the density $\rho_\kappa(x)=e^{-\kappa b(x)}$ for $0<\kappa\leq \kappa_0(\delta)$.
Let $k_\kappa$ be the quasihyperbolic metric of $X_\kappa$.

\begin{lem}\label{z14}
The identity map
$(X,d)\to (X_\kappa,k_\kappa)$
is $M$-bilipschitz with $M$ depending only on $\delta$, $K$, and $\kappa$.
\end{lem}

\bpf Let $x,y\in X$. By $(\ref{zz3})$, we get
$$k_\kappa(x,y)\leq \int_\gamma\frac{|d_{\kappa}z|}{d_\kappa (z)}=\int_0^{d(x,y)}\frac{dL_\kappa(t)}{d_\kappa(\gamma(t))}= \int_0^{d(x,y)}\frac{\rho_\kappa(\gamma(t))}{d_\kappa(\gamma(t))}\,dt\leq 2e^{10\kappa\delta}\kappa d(x,y),$$
where
$L_\kappa(t)=\int_0^t \rho_\kappa(\gamma(s))\,ds$
and $\gamma$ is a geodesic $[x,y]$ parameterized by arc length.

Thus, we only need to show the other direction. To this end, again by $(\ref{zz3})$ and Theorem \ref{q-5}, we see from \cite[(2.4)]{BHK} that
\beq\label{q0}
k_\kappa(x,y) &\geq&  \log\left(1+\frac{d_\kappa(x,y)}{d_\kappa(x)\wedge d_\kappa(y)}\right)
\\ \nonumber&\geq&  \log\left(1+\frac{\ell_\kappa([x,y])}{20C_0^2[d_\kappa(x)\wedge d_\kappa(y)]}\right)
\\ \nonumber&\geq& \log\left(1+\frac{\kappa\ell_\kappa([x,y])}{C_1[\rho_\kappa(x)\wedge \rho_\kappa(y)]}\right),
\eeq
where $C_0=e^{10\kappa\delta}$ and $C_1=C_1(K,\delta)$.

Next, we consider two cases. We first assume that $\kappa d(x,y)\leq 2\log2$. Now by (\ref{q-3}), we have
$$\rho_\kappa(z)\geq \rho_\kappa(x)e^{-\kappa d(x,z)}e^{-10\kappa\delta}\geq \frac{1}{4C_0}\rho_\kappa(x) \ \mbox{ for all $z\in[x,y]$}.$$
This leads to
$$\ell_\kappa([x,y])=\int_{[x,y]}\rho_\kappa(z)\,|dz|\geq   \frac{d(x,y)}{4C_0}\rho_\kappa(x).$$
This, together with (\ref{q0}), shows that
\beqq
k_{\kappa}(x, y)\geq \log\Big(1+\frac{\kappa d(x, y)}{4C_0 C_1}\Big)
\geq \frac{\kappa d(x, y)}{4C_0C_1+\kappa d(x, y)}
\geq \frac{\kappa d(x, y)}{4C_0C_1+2\log 2}.
\eeqq

It remains to consider the case that $\kappa d(x,y)> 2\log2$. For each $x_1\in[x,y]$, by (\ref{q-3}), we find that
\beq\label{q1}
\ell_\kappa([x,y]) &\geq&  \frac{\rho_\kappa (x_1)}{C_0}\left(\int_0^{ d(x,x_1)}e^{-\kappa t}\,dt+ \int_0^{d(x_1,y)}e^{-\kappa t}\,dt\right)
\\ \nonumber&=&  \frac{\rho_\kappa (x_1)}{C_0\kappa}\left(2-e^{-\kappa d(x,x_1)}-e^{-\kappa d(x_1,y)}\right)
\\ \nonumber&\geq&  \frac{\rho_\kappa (x_1)}{2C_0\kappa}.
\eeq

Consider the extended geodesic triangle $\Delta=\Delta(x,y,\xi)$ and take $x_1=w_y$ with $d(w_y,y)=(\xi|x)_y$. By Lemma \ref{z2}(\ref{q-1}), we get
$$b(x)\vee b(y)\geq b(w_y)+d(x,w_y)\vee d(w_y,y)-16\delta \geq b(w_y)+\frac{1}{2}d(x,y)-16\delta.$$
This implies that
\be\label{q2} \rho_\kappa(x)\wedge \rho_\kappa(y)=e^{-\kappa(b(x)\vee b(y))}\leq e^{16\kappa\delta}e^{-\frac{\kappa}{2}d(x,y)}\rho_\kappa(w_y).\ee
By (\ref{q0}), (\ref{q1}), and (\ref{q2}), it follows from \cite[(2.12)]{BHK} that
$$k_\kappa(x,y)\geq \log\left(1+\frac{t_1\rho_\kappa(w_y)}{\rho_\kappa(x)\wedge \rho_\kappa(y)}\right)\geq \log\left(1+t_2 e^{\frac{\kappa}{2}d(x,y)}\right)\geq t_2 \frac{\kappa}{2}d(x,y),$$
for some constants $t_1,t_2\in(0,1)$ depending only on $K$ and $\delta$.
\epf

\subsection{Proof of (S-\ref{qq-4})}
Let $(\Omega,d)$ be an unbounded $A$-uniform space and $k$ the quasihyperbolic metric of $\Omega$. By \cite[Proposition 3.12]{BHK}, there is a $\xi\in \partial_\infty (\Omega, k)$ such that every quasihyperbolic geodesic ray $\gamma:[0,\infty)\to \Omega$ with $\ell_d(\gamma)=\infty$ corresponding to $\gamma(\infty)=\xi$. And the natural map
$\varphi:\partial_\infty (\Omega, k)\to \partial \Omega\cup \{\infty\}$
exists which is a bijection with $\varphi(\xi)=\infty$. We denote by $(\Omega,d_\kappa)$ the conformal deformation of the space $(\Omega,k)$ associated to the density $\rho_\kappa(x)=e^{-\kappa b(x)}$, where $0<\kappa\leq \kappa_0(A)$ and $b:(\Omega,k)\to \mathbb{R}$ is a Busemann function based at $\xi$.

We need to show that the identity map $(\Omega,d)\to (\Omega,d_\kappa)$ is a quasisimilarity. To this end, we first prove a technical result.

\begin{lem}\label{z9}
Let $\gamma=[a,\xi]$ be a quasihyperbolic geodesic ray with $\gamma(0)=a\in \Omega$. For each point $v\in \Omega$ and  $u\in\gamma$ with $\ell(\gamma[a,u])=d(a,v)$, we have
\be\label{w-1} k(a,y) \doteq_{C(A)} (v|\xi)_a.\ee
\end{lem}
\bpf For an integer $n$ sufficiently large we may assume that
$d(a,\gamma(n))\geq 2d(a,v).$
Put $z_n=\gamma(n)$. Take a quasihyperbolic geodesic $[a,v]$ connecting $a$ to $v$ in $\Omega$. By a direct application of \cite[Lemma 2.30]{Vai} in the setting of uniform spaces, we know that
\be\label{z-0-1} k(u,z_n) \doteq_{C(A)} \dist_k(z_n,[a,v]),\ee
where $\dist_k(z_n,[a,v])$ is the quasihyperbolic distance of $z_n$ to $[a,v]$.
Because $\Omega$ is $A$-uniform, we know from
Theorem \ref{thm-3} that $(\Omega,k)$ is $\delta$-hyperbolic with $\delta=\delta(A)$. By \cite[(3.2)]{BHK}, we get
\be\label{z-0-2} \dist_k(z_n,[a,v]) \doteq_{C(A)} (a|v)_{z_n}.\ee
Therefore, by (\ref{z-0-1}) and (\ref{z-0-2}), we obtain
\be\label{b-5} (v|z_n)_a=k(a,z_n)-(v|a)_{z_n}  \doteq_{C(A)} k(a,z_n)-k(u,z_n)=k(a,u).\ee
Because $\{\gamma(n)\}\in\xi$, by (\ref{b-5}) and Lemma \ref{z0}, we obtain (\ref{w-1}).
\epf

\begin{lem}\label{z10}
The identity map $(\Omega,d)\to (\Omega,d_\kappa)$ is quasisimilar with the parameters $\theta, \lambda$, and $\tau$ depending only on $A$.
\end{lem}
\bpf It follows from Theorem \ref{thm-3} that $(\Omega,k)$ is $\delta$-hyperbolic for some $\delta=\delta(A)\geq 0$. We shall show that there is a number $\tau=\tau(A)\in(0,1/3)$  such that for all $x\in\Omega$ and for any $y,z\in B(x,\tau d(x))$,
\be\label{s-2} \kappa_0 k(y,z)\leq 1\;\;\mbox{and}\;\; (y|z)_b \doteq_{C(A)} b(x). \ee
Choose $t_0\in(0,1)$ with
$$e^{t_0/{\kappa}}-1=\tau/2.$$
Because $y,z\in B(x,\tau d(x))$, we obtain by the triangle inequality that
$$d(y,z)\leq 2\tau d(x)$$ and $$d(y)\wedge d(z)\geq (1-\tau)d(x).$$
As $(\Omega,d)$ is $A$-uniform, we see from \cite[(2.16)]{BHK} that
$$k(y,z)\leq 4A^2\log\left(1+ \frac{d(y,z)}{d(y)\wedge d(z)}\right)\leq 4A^2\log\left(1+\frac{2\tau}{1-\tau}\right).$$
It follows that there is a constant $\tau=\tau(A)\in(0,1)$ such that $\kappa_0 k(y,z)\leq 1$ and the first assertion in (\ref{s-2}) is true.

Now we prove the second assertion in (\ref{s-2}). This can be seen from (\ref{b-1}) and the first assertion of (\ref{s-2}). Indeed, on the one hand, we have
$$(y|z)_b\leq b(y)+5\delta\leq b(x)+k(x,y)+15\delta\leq b(x)+1/{\kappa_0}+15\delta$$
and, on the other hand, a similar argument as above shows that
\beq\nonumber
(y|z)_b  & \geq & b(y)-k(y,z)-5\delta
\\ \nonumber & \geq&  b(x)-k(x,y)-k(y,z)-15\delta
\\ \nonumber & \geq & b(x)-2/{\kappa_0}-15\delta.
\eeq
This proves (\ref{s-2}).

By (\ref{s-2}) and (\ref{zz2}), it follows that for all $x\in\Omega$ and $y,z\in B(x,\tau d(x))$,
\be\label{s0} d_\kappa (y,z) \asymp_{C(A)}e^{-\kappa (y|z)_b}k(y,z) \asymp_{C(A)} \rho_\kappa(x) k(y,z).\ee
We are ready to check the second requirement (QS-2) of quasisimilarity with the data $(\lambda, \tau)$ depending only on $A$. By \eqref{s0},
we only need to show that
\be\label{s0a}k(y,z)\asymp_{C(A)}\frac{d(y,z)}{d(x)} \ \mbox{ for all $x\in\Omega$ and $y,z\in B(x,\tau d(x))$,}\ee
where $\tau$ is the constant such that (\ref{s-2}) holds. We observe that
$d(y)\wedge d(z)\geq (1-\tau)d(x)$
and
$$d(y,z)\leq 2\tau d(x)\leq \frac{2\tau}{1-\tau}d(y)\leq d(y).$$
Again by \cite[(2.16) and (2.4)]{BHK} and the above two inequalities, it follows easily that
$$k(y,z)\leq 4A^2\log\left(1+ \frac{d(y,z)}{d(y)\wedge d(z)}\right)\leq  \frac{4A^2 d(y,z)}{(1-\tau)d(x)}$$
and
$$k(y,z)\geq \log\left(1+ \frac{d(y,z)}{d(y)\wedge d(z)}\right)\geq   \frac{d(y,z)}{d(y,z)+d(y)} \geq\frac{d(y,z)}{(1+3\tau)d(x)}.$$
Therefore, the second requirement (QS-2) of quasisimilarity  for the identity map $(\Omega,d)\to (\Omega,d_\kappa)$ is satisfied.

It remains to check the $\theta$-quasisymmetry of the identity map $\Omega\to \Omega_\kappa$ with $\theta$ depending only on $A$. By \cite[Theorem $6.6$]{Vai-4}, it suffices to show that for each triple of distinct points $x,y,z\in \Omega$ the following statement holds:
\be\label{s1} d(x,y)\leq d(x,z) \;\;\;\;\;\mbox{implies}\;\;\;\;\; d_{\kappa}(x,y) \leq H d_{\kappa}(x,z),\ee
because both $\Omega$ and $\Omega_\kappa$  are uniform and thus quasiconvex.

By (\ref{zz2}), we have
\be\label{s-3} \frac{d_{\kappa}(x,y)}{d_{\kappa}(x,z)} \asymp_{C_1(A)} e^{\kappa[(x|z)_b-(x|y)_b]} \frac{1\wedge (\kappa k(x,y))}{1\wedge (\kappa k(x,z))}.\ee
Take a quasihyperbolic geodesic ray $[x,\xi]$ connecting $x$ to $\xi$ in $\Omega$. For each $v\in\{y,z\}$, Lemma \ref{z9} ensures that there is a point $u_v$ in $[x,\xi]$ corresponding to the point $v$ such that
$$d(x,v)=\ell(\gamma[x,u_v]) \;\;\;\;\;\mbox{and}\;\;\;\;\; k(x,u_v) \doteq_{C_2(A)} (v|\xi)_x.$$
Because $d(x,y)\leq d(x,z)$, we get
$$(z|\xi)_x-(y|\xi)_x \doteq_{C_3(A)} k(x,u_z)-k(x,u_y)=k(u_y,u_z).$$
This, together with (\ref{zz0.1}) and Lemma \ref{z0}, implies that
\beq\label{s-4}
(x|z)_b-(x|y)_b &\doteq_{C_4(A)}& (x|z)_o-(z|\xi)_o-(x|y)_o+(y|\xi)_o
\\ \nonumber&\doteq_{C_5(A)}& (y|\xi)_x-(z|\xi)_x
\\ \nonumber&\doteq_{C_6(A)}& -k(u_y,u_z)\leq 0.
\eeq
Furthermore, it follows  from (\ref{s-3}) and (\ref{s-4}) that there is a constant $C_7=C_7(A)\geq 1$ such that
$$\frac{d_{\kappa}(x,y)}{d_{\kappa}(x,z)} \leq C_7 \frac{1\wedge (\kappa k(x,y))}{1\wedge (\kappa k(x,z))}.$$

If $\kappa k(x,z)>t_0$, by the above inequality, the statement (\ref{s1}) follows. If $\kappa k(x,z)\leq t_0$, then we obtain by our choice of $t_0$ that
$$d(x,y)\leq d(x,z)\leq (e^{k(x,z)}-1)d(x)\leq \frac{\tau}{2} d(x).$$
Thus we have $y,z\in B(x,\tau d(x))$. The desired inequality (\ref{s1}) follows from (\ref{s0}) and \eqref{s0a} in this case.
\epf

The proof of Theorem \ref{main thm-1} is complete.
\qed
\bigskip
%


\end{document}